\documentclass  [a4paper, 12pt] {article}
\usepackage{amsmath}
\usepackage{amssymb}
\usepackage{amsfonts}
\usepackage [latin1] {inputenc}
\usepackage [T1] {fontenc}
\usepackage[french]{babel}

\def\R{\mathbb{R}}

\def\N{\mathbb{N}}
\def\Z{\mathbb{Z}}

\title{ SEMICIRCLE LAW FOR RANDOM MATRICES OF LONG-RANGE PERCOLATION MODEL  }

\author{{\bf S.Ayadi}\thanks{LMV - Laboratoire de Math\'ematiques
de Versailles, Universit\'e de Versailles Saint-Quentin-en-Yvelines,
78035 Versailles (FRANCE). E-mail: ayadi@math.uvsq.fr}}
\date{}
\begin{document}
\maketitle

{\bf Abstract}: We study the normalized eigenvalue counting measure $d\sigma$ of matrices of long-range percolation model. These are $(2n+1)\times (2n+1)$ random real symmetric matrices $H=\{H(i,j)\}_{i,j}$ whose elements are independent random
variables taking zero value with probability $1-\psi\left( (i-j)/b\right)$,
$b\in \mathbb{R}^{+}$, where $\psi$ is an even positive function
$\psi(t)\le{1}$ vanishing at infinity. It is shown that if the third moment of $\sqrt{b}H(i,j), \ i\le{j}$ is uniformly bounded then the measure $d\sigma:=d\sigma_{n,b}$ weakly converges in probability in the limit $n,b\rightarrow\infty$, $b=o(n)$ to the semicircle (or Wigner) distribution. The proof uses the resolvent technique combined with the cumulant expansions method. We show that the normalized trace of resolvent $g_{n,b}(z)$ converges in average and that the variance of $g_{n,b}(z)$ vanishes. In the second part of the paper, we estimate the rate of decreasing of the variance of $g_{n,b}(z)$, under further conditions on the moments of $\sqrt{b}H(i,j), \ i\le{j}$.

\vskip0,5cm

 \vskip0,2cm

 \vskip0,2cm

{\bf running title}: {\it Semicircle law for percolation model.}

\section{Introduction} Spectral theory of random matrices is a relatively new branch of mathematics that intensively develops due to its rich mathematical content and aho due to numerous applications in theoretical phisics, wireless communications, financial mathematics and other fields (see reviews monographs \cite{AAA}, \cite{CC}, \cite{HHHHH}, \cite{HHHH} and references therein).

In theoretical physics random matrices of infinitely increasing dimensions are used
to give the statistical description of systems with large number of degrees of freedom.

 The first application was in nuclear physics, where E.Wigner
proposed to consider the eigenvalues of real symmetric random
matrices as a model for the energy levels of heavy atomic nuclei
(see the monograph \cite{CC}). The
real symmetric random matrix $A_{n}$ of size $2n+1$ introduced by E. Wigner is
defined by

$$
 A_{n}(i,j)=N^{-1/2}a(i,j), \quad |i|, |j|\le{n}, \eqno (1.1)
$$
where $N=2n+1$ and $\{a(i,j), \ i\le{j}\}$ are independent and identically distributed random variables defined
on the same probability space $(\Omega,\mathfrak{F},\mathbb{P})$
such that
$$
{\bf E}\{a(i,j)\}=0, \quad {\bf E}\{a^{2}(i,j)\}:=v^{2}, \eqno (1.2)
$$
where ${\bf E}$ is the mathematical expectation with respect to $\mathbb{P}$.

Denoting by $\lambda^{(n)}_{-n}\le{...}\le{\lambda^{(n)}_{n}}$ the
eigenvalues of $A_{n}$, the normalized eigenvalue counting function
is defined by
$$
\sigma_{n}(\lambda,A_{n}):=N^{-1}\sharp\{\lambda^{(n)}_{j}\le{\lambda}\}.
\eqno (1.3)
$$
E.Wigner \cite{BB} proved that in the case when $a(i,j)$ has all
moments finite, the eigenvalue counting measure
$d\sigma_{n}(\lambda,A_{n})$ weakly converges in average as
$n\rightarrow \infty$ to a distribution
$d\sigma_{sc}(\lambda)$, where the nondecreasing function
$\sigma_{sc}(\lambda)$ is differentiable and its derivative
$\rho_{sc}$ writes as follows
$$
\rho_{sc}{(\lambda)}= \sigma_{sc}^{'}(\lambda) =  \left\{
\begin{array}{lll}
\frac{\sqrt{4v^{2}-\lambda^{2}}}{2\pi v^{2}} & \textrm{if} &
|\lambda|\le{2v} \\
0                               &  \textrm{ }  &    otherwise  .
\end{array}\right.
\eqno (1.4)
$$

This limiting distribution $(1.4)$ is referred to as the Wigner
distribution, or the semicircle law.

  Since then, the convergence to the semicircle distribution was proved for
various random matrix ensembles that generalize the Wigner ensemble
\cite{I,E,H}. Among them we cite the band and the dilute random
matrix ensembles \cite{F,DDD,CCC,AA,BBB}.

 In the band random matrices model, the matrix elements take zero
value outside of the band of width $b_{n}$ along the principal diagonal, for some positive sequence of real numbers. This ensemble can be obtained from $A_{n}$ $(1.1)$ by multiplying each $a(i,j)$ by $I_{(-1/2,1/2)}\left( (i-j)/b_{n}\right)$, where
$$
I_{A}(t) = \left\{
\begin{array}{lll}
1 & \textrm{if} & t\in{A}, \\
0 &  \textrm{if} & t\in\mathbb{R}\setminus A
\end{array}\right.
$$
is the indicator function of the interval $A$.
The ensemble of dilute random matrices can be obtained from $A_{n}$ $(1.1)$ by multiplying $a(i,j)$ by independent Bernoulli random variables 
$$
\widehat{d}(i,j) = \left\{
\begin{array}{lll}
1 & \textrm{with probability} & p_{n}/n, \\
0 &  \textrm{with probability} & 1 - p_{n}/n, \quad 0\le{p_{n}}\le{n}.
\end{array}\right. 
$$
Assuming that $b(n)=o(n)$ for large $n$, the semicircle law is observed, after proper renormalization,for both ensembles, in the limit $b_{n}\rightarrow\infty$ \cite{AA} and
$p_{n}\rightarrow\infty$ \cite{DDD} as $n\rightarrow\infty$.

 It should be noted that the matrix $\{\widehat{d}(i,j)\}_{i,j}$ represents the adjacency matrix of random Erd\H os-Rényi graphs (see \cite{BBBB} for definitions and more details). Random Erd\H os-Rényi graphs are extensively studied. Recent applications also involve the problems of spectral theory of random graphs models \cite{EEE} and references therein.

 Our work is motivated by papers where the random graphs of certain percolation models are considered. More precisely,we are interested in a generalization of the tow ensembles mentioned above. Roughly speaking, we consider band random
matrices with a random width.
  To proceed, we multiply each matrix element $a(i,j)$ by $b_{n}^
{-1/2}d(i,j)$, where $0<b_{n}\le{N}$ is a sequence of integers and
$d(i,j), \ |i|,\ |j|\le{n}$ is a Bernoulli random variable with
$$
d(i,j) = \left\{
\begin{array}{lll}
1 & \textrm{with probability} & \psi{(\frac{i-j}{b_{n}})}, \\
0 &  \textrm{with probability} & 1 - \psi{(\frac{i-j}{b_{n}})}
\end{array}\right.  \eqno (1.5)
$$
and $0\le{\psi(t)}\le{1}$ is an even function vanishing as
$t\rightarrow\infty$. We assume that $\mathcal{D}_{n}=\{d(i,j),
\quad |i|,\ |j|\le{n}\}$ is a family of independent random
variables, also independent from ${\cal A}_{n}= \{a(i,j), \quad
|i|,\ |j|\le{n}\}$.

 The family $\{d(i,j), \quad |i|,\ |j|\le{n}\}$ can be regarded as the
adjacency matrix of the family of random graphs $\{\Gamma_{n}\}$
with $N=2n+1$ vertices $(i,j)$ such that the average number of edges
attached to one vertex is $b_{n}$. Hence, each edge $e(i,j)$ of the
graph is present with probability $\psi\left( (i-j)/b\right)$ and not
present with probability $1-\psi\left( (i-j)/b\right)$. Below are some well
known examples:

\begin{enumerate}

\item[$-$] In theoretical physics, the ensemble $\{\Gamma_{n}\}$ with $\psi(t)=e^{-|t|^{s}}$
is referred to as the Long-Range Percolation Model (see for exemple \cite{K} and references therein). Our
ensemble can be regarded as a modification of the adjacency matrices
of $\{\Gamma_{n}\}$. To our best knowledge, the spectral properties of
this model has not been studied yet.

 \vskip0,2cm

\item[$-$] It is easy to see that if one takes $b_{n}=N$ and $\psi\equiv 1$,
then one recovers the Wigner ensemble $(1.1)$.

\vskip0,2cm
\item[$-$]  If one considers $(1.5)$ with $\psi(t)=I_{(-1/2,1/2)}(t)$, one gets the band random matrix ensemble described before
\cite{AA}.

\end{enumerate}

 The paper is organized as follows. In section 2, we define the
random matrix ensemble $H_{n}$ of long-range percolation model and
state our main results. In section 3, we use the cumulant expansions
method to derive relations for ${\bf E}\{g_{n}(z)\}$, where
$g_{n}(z)$ is the Stieltjes transform of the normalized eigenvalue
counting measure of our random matrix ensemble. In section 4, we
show that ${\bf Var}\{g_{n}(z)\}$ vanishes and prove the semicircle
law for the random matrix $H_{n}$ under the condition that
$\sup_{i,j} {\bf E}\{|a(i,j)|^{3}\}<\infty$. In section 5, assuming
more conditions on the moments of $a(i,j)$, we derive more precise
and powerful estimates for the variance ${\bf Var}\{g_{n}(z)\}$.
\vskip0,5cm

\section{The ensemble and main results}
 Let us consider a family of independent real random variables
${\cal A}_{n}= \{a(i,j):  |i|, |j|\le{n} \}$ satisfying $(1.2)$ and
defined on the same probability space
$(\Omega,\mathbb{F},\mathbb{P})$.

 Let $\psi(t)$, $t\in\mathbf{R}$, be a real continuous even function such that:

$$
0\le{\psi{(t)}}\le{1}    \    and \
       \int_{\mathbb{R}}\psi{(t)}dt = 1.
\eqno (2.1)
$$

Introduce a family of independent Bernoulli random variables
 ${\cal D}_{n}=\{d(i,j):  |i|, |j|\le{n} \}$ as in $(1.5)$, that is
independent of the family ${\cal A}_{n}$. We assume that ${\cal
A}_{n}$ and ${\cal D}_{n}$ are defined on the same probability
space. \vskip0,2cm
 Define a real symmetric $N\times N$ random matrix $H_{n}$ by:
$$
H_{n}(i,j)=b_{n}^{-1/2}a(i,j)d(i,j), \quad i\le{j}, \quad |i|,
|j|\le{n}, \eqno (2.2)
$$
where $0<b_{n}\le{N}$, $N=2n+1$. We study asymptotic spectral
properties of the ensemble $\{H_{n}\}$ in the limit when
$n\rightarrow\infty$ with
$$
\lim_{n\rightarrow\infty} b_{n}=\infty \quad and \quad
\lim_{n\rightarrow\infty} \frac{b_{n}}{n}=0. \eqno (2.3)
$$
This limit $(2.3)$ when $1<<b_{n}<<n$ is known as the limit of
relativity narrow band width.

 Here and below the family $\{H_{n}\}$ is referred to as the ensemble of random matrices of long-range percolation model.

Let us introduce
$$
\hat{\mu}_{r}=\sup_{i,j}{\bf E}\{|a(i,j)|^{r}\},  \  r\in
\N, \eqno (2.4)
$$
the absolute moments of $a(i,j)$.
\vskip0,5cm
Our main result is as follows.

\vskip0,5cm

{\bf Theorem 2.1 } {\it If $\hat{\mu}_{3}<\infty$ $(2.4)$, the normalized eigenvalue counting function
$\sigma_{n}(\lambda,H_{n})$ in $(1.3)$ for the ensemble given in
$(2.2)$ with $\psi$ satisfies $(2.1)$, converges in probability in
the limit $(2.3)$ to a nonrandom distribution function:
$$
p-\lim _{n\rightarrow\infty} \
\sigma_{n}{(\lambda,H_{n})}=\sigma_{sc}(\lambda), \eqno (2.5)
$$
where $\sigma_{sc}(\lambda)$ is given by $(1.4)$.}

\vskip0,5cm

Let $\{\mu_{n}\}_{n}$ be a sequence of probability measures that
converge weakly to a probability measure $\mu$. Let $g_{n}$ and $g$
denote the Stieltjes transforms of $\mu_{n}$ and $\mu$ respectively
defined by
$$
g_{n}(z) = \int \frac{\mu_{n}(d\lambda)}{\lambda-z}; \ g(z)=\int
\frac{\mu(d\lambda)}{\lambda-z}, \eqno (2.6)
$$
for $z\in \mathbb{C}\setminus \mathbb{R}$. Then the weak convergence of
$\mu_{n}$ is equivalent to the convergence of the Stieltjes
transform (in the uniform topology on compact sets \cite{HHH}).

Let $\sigma_{n}(d\lambda,H_{n})$ be counting measure of the
eigenvalue of $H_{n}$:
$$
\sigma_{n}(d\lambda,H_{n})=\frac{1}{N}\sum_{|i|\le{n}}\delta_
{\lambda_{i}}(d\lambda), \eqno (2.7)
$$
where $\delta_ {\lambda_{i}}$ is the Dirac mass at $\lambda_{i}$.
Then Theorem 2.1 is equivalent to the following statement.

\vskip0,5cm

 {\bf Theorem 2.2} {\it Let $g_{n}(z)$ be the Stieltjes transform of
the normalized eigenvalue counting measure
$\sigma_{n}(d\lambda,H_{n})$ $(1.3)$. If $
\hat{\mu}_{3}<\infty$ $(2.4)$, then in the limit $(2.3)$, there
exists a nonrandom function $w(z)$ such that for all nonreal $z$

$$
p-lim _{n\rightarrow \infty} \ g_{n}(z) = w_{sc}(z), \eqno (2.8)
$$

where $w_{sc}(z)$ is the Stieltjes transform of the semicircle law
given by $(2.5)$. }

\vskip0,5cm

{\bf Remarks} \vskip0,1cm
\begin{enumerate}

\item The limiting function $w_{sc}(z)$ is the solution of equation
$$
w_{sc}(z)=\frac{1}{-z-v^{2} w_{sc}(z)} \eqno (2.9)
$$
that maps $\mathbb{C}^{+}\rightarrow \mathbb{C}^{+}$ when
$$
\mathbb{C}^{+}=\{z\in \mathbb{C}: \quad Im z>0\}
$$
and the parameter $v$ is determined by $(1.2)$.
\item To prove $(2.8)$, we start by proving it for $z\in
\Lambda_{\eta}$, where

$$
 \Lambda_{\eta} = \{z\in\mathbf{ C}:
 \  \eta \le{ |Im{(z)}|}\} \quad with \ \eta=2v+1.
\eqno (2.10)
$$
\end{enumerate}

The convergence in probability $(2.8)$ is shown in a standard way:

first we prove that ${\bf E}\{g_{n}(z)\}$ converges to $w_{sc}(z)$ $(2.9)$
and that ${\bf Var}\{g_{n}(z)\}$ vanisfes. More precisely, we are going to prove the following statement.

\vskip0,5cm

{\bf Theorem 2.3} {\it If  $\hat{\mu}_{3}<\infty$ $(2.4)$ and
$z\in\Lambda_{\eta}$ $(2.10)$, then in the limit
$n\rightarrow\infty$ $(2.3)$,
$$
\lim_{n\rightarrow\infty}{\bf E}\{g_{n}(z)\}=w_{sc}(z) \eqno (2.11)
$$
and
$$
{\bf Var}\{g_{n}(z)\} = O(b_{n}^{-1/2}), \eqno (2.12)
$$
where $w_{sc}(z)$ is given by equation $(2.9)$.

 }
\vskip0,5cm
Theorem 2.2 follows then using Tchebychev's Inequality
and gives Theorem 2.1.

\vskip1cm
   To prove Theorem 2.3, we mainly follow papers \cite{A,B}. We apply the
cumulant expansions method of \cite{B} and get similar relations to
the ones obtained in \cite{A}. So, combining these two approaches,
we develop a method to study the resolvent of $H_{n}$. It is rather general and can be used to study other random matrix ensemble with jointly independent entries.

  It is not hard to show that the estimate $(2.11)$ together with
$(2.12)$ imply convergence in probability $(2.8)$. We give the
details at the end of section 3. Then $(2.8)$ implies convergence
$(2.5)$ ( see section 3  for the proof).

 Then in section 3, we give the limiting equation from
${\bf E}\{g_{n}(z)\}$. In section 4, we prove Theorem 2.3 and we
deduce from this the results of Theorem 2.2 and Theorem 2.1. In
section 5, we give more details on asymptotic behavior of the
variance and prove auxiliary statements.

 For the sake of simplicity we omit the subscript $n$ in $b_{n}$, $G_{n}$ and
$H_{n}$, so that $G_{n}=G$ $(2.12)$ is the resolvent of the matrix
$H_{n}=H$ $(2.2)$.

\vskip0,5cm

\section{Limiting equation for ${\bf E}\{g_{n}(z)\}$}
The resolvent
$$
G_{n}(z)=(H_{n}-zI)^{-1}, \quad Imz\neq{0},
$$
is widely exploited in the spectral theory of operators. Its normalized trace coincides with the Stieltjes transform $g_{n}(z)$ $(2.6)$ of
the normalized eigenvalue counting measure $\sigma_{n}(d\lambda,H_{n})$ $(1.3)$;
$$
g_{n}(z)=\int \frac{\sigma_{n}(d\lambda,H_{n})}{\lambda-z}=\frac{1}{N}\sum_{|j|\le{n}}\frac{1}{\lambda_{j}^{(n)}-z}.
$$

Let us consider the resolvent identity for two hermitian matrices
$H$ and $H'$:
$$
G - G' = G(H - H')G', \eqno (3.1)
$$
where $G(z) = (H - zI)^{-1}$,  $G^{'}(z) = (H' - zI)^{-1}$.

Regarding $(3.1)$ with $H=H_{n}$ $(2.2)$ and $H' = 0$, we obtain
relation
$$
G(i,j) = \xi\delta_{ij} - \xi\sum_{|p|\le{n}} G(i,p)H(p,j), \quad
\xi\equiv{-z^{-1}},  \eqno (3.2)
$$
where $\delta$ denotes the Kronecker symbol

$$
\delta_{ij} = \left\{
\begin{array}{lll}
1 & \textrm{if} & i=j \\
0 &  \textrm{if} & i\neq{j},
\end{array}\right.
$$
$H(i,j)=H_{n}(i,j), \ |i|, |j|\le{n}$ are the entries of the matrix
$H_{n}$ and $G(i,j)$ are the entries of the resolvent $G_{n}=(H_{n}-zI)^{-1}$.

 Relation $(3.2)$ implies that
$$
{\bf E}G(i,i) = \xi - \xi\sum_{|p|\le{n}} {\bf E}\{G(i,p)H(p,i)\}.
\eqno (3.3)
$$
To compute the mathematical expectation ${\bf E}\{G(i,p)H(p,i)\}$,
we use the cumulants expansion method proposed in \cite{B}.

\subsection{Cumulant expansions and principal equation}
Let us start with the description of our basic technical tools.

  $(i)$. {\it The cumulant expansions formula}. Let us consider a
family $\{X_{j}, \ j=1,..,m\}$ of independent real random variables
defined on the same probability space such that ${\bf
E}\{|X_{j}|^{q+2}\}<\infty $ for same $q\in \N$ and
$j=1,..,m$. Then for any complex-valued function $F(t_{1},..,t_{m})$
of the class $\mathcal{C}(\mathbf{R}^{m})$ and for all $j$, one has
$$
{\bf E}[X_{j} F(X_{1},..,X_{m})] = \sum_{r=0}^{q} \frac{K_{r+1}}{r!}
{\bf E}\{\frac{\partial^{r}F(X_{1},..,X_{m})}{\partial{X_{j}^{r}}}\}
  +  \epsilon^{(j)}_{q},
\eqno (3.4)
$$
where $K_{r}$ is the r-th cumulant of $X_{j}$ \cite{B} and the
remainder $\epsilon_{q}$ can be estimated by the inequality
$$
|\epsilon^{(j)}_{q}|\le{C \sup_{t\in\mathbf{R}^{m}}
|\frac{\partial^{q+1}F(t)}{\partial{t_{j}^{q+1}}}| {\bf
E}\{|X_{j}|^{q+2}\}} \eqno (3.5)
$$
and $C$ is a constant depending only on $q$.

Relations $(3.4)$ and $(3.5)$ can be proved by using the Taylor's formula
(see e.g.\cite{B} for example). At the end of section 5, we prove $(3.4)$ in the case of $q=5$ and give more details on the form of the remainder terms $\epsilon_{q}^{j}$ of $(3.4)$.
\vskip0,2cm
 It is known that the cumulants can be expressed in terms of the moments of $X_{j}$. If we denote ${\bf E}(X_{j})=
{\bf E}(X_{j}^{3})={\bf E}(X_{j}^{5})=0$ and $\mu_{r}={\bf
E}(X_{j}^{r})$ with $j=1,..,m$, one obtains: $K_{1}=\mu_{1}=0$,
$K_{2}=\mu_{2}$, $K_{3}=0$, $K_{4}=\mu_{4}-3\mu^{2}_{2}$, $K_{5}=0$,
$K_{6}=\mu_{6}-15\mu_{4}\mu_{2}+30\mu^{3}_{2}$, etc.

\vskip0,1cm

 $(ii)$. If $H$ is a real symmetric $N\times N$ matrix with elements $H(i,j),
\ |i|, |j|\le{n}$, and $G=(H-z)^{-1}(i,j)$ is its resolvent, then

$$
\frac{\partial{G_{st}}}{\partial{H_{jk}}} = \left\{
\begin{array}{lll}
- G_{sj}G_{kt} & \textrm{si} & j=k \\
- G_{sj}G_{kt} - G_{sk}G_{jt} &  \textrm{si} & j\neq{k}.
\end{array}\right.
\eqno (3.6)
$$
 The relation $(3.6)$ follows from the resolvent identity $(3.1)$.
\vskip0,5cm

Now, let us explain the main idea of the proof of relation $(2.15)$.
Applying $(3.4)$ to the resolvent identity $(3.3)$, one gets
$$
{\bf E}\{G(i,p)H(p,i)\}=\frac{v^{2}}{b}{\bf
E}\{\frac{\partial{G(i,p)}}{\partial{H(p,i)}}
\}\psi(\frac{p-i}{b})+\epsilon_{ip}.
$$
Substituting this equality in $(3.3)$ and using $(3.6)$, we arrive at relation 
$$
{\bf E}G(i,i)=\xi+\frac{\xi v^{2}}{b}\sum_{|p|\le{n}}{\bf
E}\{G(i,i)G(p,p)\}\psi(\frac{i-p}{b})+R_{n,b}(i),
$$
where $R_{n,b}(i)$ vanishes as $n,b\rightarrow\infty$ $(2.3)$ (see
subsection 3.2 for more details).

 If one assumes that the average ${\bf
E}\{G(i,i)b^{-1}\sum_{p}G(p,p)\psi([p-i]/b)\}$ factorizes, then
$g_{n,b}(i)={\bf E}G(i,i)$ satisfies equality
$$
g_{n,b}(i)=\xi+\frac{\xi
v^{2}}{b}g_{n,b}(i)\sum_{|p|\le{n}}g_{n,b}(p)\psi(\frac{p-i}{b})+R_{n,b}(i).
$$
 Assuming that the value
$b^{-1}\sum_{p}\psi([i-p]/b)g_{n,b}(p)$ does not depend on $i$ in
the limit $n,b\rightarrow\infty$, we get convergence
$g_{n,b}(i)\rightarrow w(z)$, for all $|i|\le{n-bL}$, where $L$ is
sufficiently large and $w(z)$ is the solution of equation (c.f.$(2.9)$)
$$
w(z)=\xi + \xi v^{2}w(z)^{2}, \quad \xi\equiv{-z^{-1}}.
$$

\subsection{Derivation of relations for ${\bf E}G(i,i)$}
Let us consider the average ${\bf E}\{G(i,p)H(p,i)\}$. For each pair
$(i,p)$, $G(i,p)$ is a smooth function of $H(p,i)$. Its
derivatives are bounded because of equation $(3.6)$ and the
inequality

$$
|G(i,p)|\le{||G||}\le{\frac{1}{|Imz|}}, \eqno (3.7)
$$
which holds for the resolvent of any real symmetric matrix. Here and below we denote by $||e||_{2}^{2}=\sum_{|i|\le{n}}|e(i)|^{2}$ the euclidean norm of
the vector $\{e(i)\}_{|i|\le{n}}$ and $
||G||=\sup_{||e||_{2}=1}||Ge||_{2}$.

  Inequality $(3.7)$ implies that, $|D^{2}_{pi}G_{ip}| \leq{M|Im z|^{-3}}$
where $M$ is an absolute constant. In what follows, we use the notation
$D^{r}_{pi}=\partial^{r}/\partial{H^{r}(p,i)}$.

  According to $(2.2)$, $(1.5)$ and the condition $\hat{\mu}_{3}<\infty$,
the third absolute moment of $H(p,i)$ is of order $b^{-3/2}$. Applying $(3.4)$ to ${\bf E}\{G(i,p)H(p,i)\}$ with $q=1$, we get relation

$$
{\bf E}\{G(i,p)H(p,i)\} =K_{1}{\bf E}( G(i,p))+K_{2}{\bf
E}(D^{1}_{pi}G(i,p))+\epsilon _{ip} \eqno (3.8)
$$
where $K_{r}$ is the r-th cumulant of $H(p,i)$ and $\epsilon_{ip}$
satisfies the estimate

$$
|\epsilon_{ip}|\le{\frac{C}{b^{3/2}}\sup_{i,p}|D^{2}_{pi}G(i,p)|{\bf
E}\{|a(p,i)|^{3}\}{\bf E}\{|d(p,i)|^{3}\}}
$$
$$
\le{\frac{CM}{b^{3/2}|Imz|^{3}}{\bf E}\{|a(p,i)|^{3}\}{\bf
E}\{|d(p,i)|^{3}\}}, \eqno (3.9)
$$
for some constants $C$ and $M$.

To compute the partial derivatives, we use $(3.7)$. For a complex symmetric matrix $G$, one has
$D^{a}_{ii}G(i,i)=(-1)^{a}a!G(i,i)^{a+1}$, with $a=1,2$ and
$$
D^{1}_{pi}G(i,p) = - G(i,p)^{2} - G(i,i)G(p,p), \eqno (3.10)
$$

$$
D^{2}_{pi}G(i,p) = 2G(i,p)^{3} + 6G(i,i)G(p,p)G(i,p), \eqno (3.11)
$$
for distinct $i$ and $p$. Using formulas $(3.4)$ and $(3.5)$, the
cumulants can be expressed in terms of the moments. Considering
$X=H(p,i)$, we get
$$ k_{1}=
\mu_{1} = {\bf E}(X) = \frac{1}{\sqrt{b}}{\bf E}\{a(p,i)\}{\bf
E}\{d(p,i)\} = 0, \eqno (3.12)
$$
$$
K_{2} = \mu_{2} = {\bf E}(X^{2}) = \frac{v^{2}}{b}{\bf
E}\{d(p,i)^{2}\}=\frac{v^{2}}{b}\psi(\frac{p-i}{b}). \eqno (3.13)
$$
Using $(3.10)$, $(3.12)$ and $(3.13)$, we rewrite $(3.8)$ in the
form
$$
{\bf E}\{G(i,p)H(p,i)\}=-\frac{v^{2}}{b}{\bf E}\{G(i,p)^{2}+
G(i,i)G(p,p)\}\psi(\frac{p-i}{b})+\epsilon_{ip}, \ p\neq i \eqno
(3.14a)
$$
$$
{\bf E}\{G(i,i)H(i,i)\}=-\frac{v^{2}}{b}{\bf
E}\{G(i,i)^{2}\}\psi(0)+\epsilon_{ii}. \eqno (3.14b)
$$
Substituting $(3.14)$ in $(3.3)$, we obtain equality
$$
{\bf E}G(i,i) = \xi +\frac{\xi v^{2}}{b} \sum_{|p|\le{n},\ p\neq i}
{\bf E}\{G(i,p)^{2}\}\psi(\frac{p-i}{b})
$$

$$
+ \frac{\xi v^{2}}{b}  \sum_{|p|\le{n},\ p\neq i}{\bf
E}\{G(i,i)G(p,p)\} \psi(\frac{p-i}{b}) +\frac{\xi v^{2}}{b}{\bf
E}\{G(i,i)^{2}\}\psi(0)- \xi\sum_{|p|\le{n}}\epsilon_{ip}. \eqno
(3.15)
$$
Relations $(3.9)$ and $(3.11)$ imply the following estimate
$$
|\epsilon_{ip}|\le{\frac{8C\hat{\mu}_{3}}{|Imz|^{3}b^{3/2}}
\psi(\frac{p-i}{b})} \eqno (3.16)
$$
for some constant $C$.

For a given a random variable, we write $\zeta^{0}$ for:
$$
\zeta^{0}=\zeta-{\bf E}\zeta,
$$
so that $(3.15)$ takes the form
$$
{\bf E}G(i,i) = \xi + \frac{\xi v^{2}}{b}{\bf E}G(i,i)
\sum_{|p|\le{n}}{\bf E}G(p,p)\psi{(\frac{i-p}{b})}
$$
$$
+ \frac{1}{b} \phi_{1}{(i)} + \phi_{2}{(i)} - \xi\sum_{|p|\le{n}}
\epsilon_{ip}, \eqno (3.17)
$$
where
$$
 \phi_{1}{(i)} = \xi v^{2}\sum_{|p|\le{n}}{\bf E}\{G(i,p)^{2}\}
 \psi{(\frac{i-p}{b})}-\frac{\xi v^{2}}{b}{\bf E}\{G(i,i)^{2}\}\psi(0),
  \quad |i|\le{n},
\eqno (3.18)
$$
$$
\phi_{2}{(i)} =\frac{\xi v^{2}}{b}\sum_{|p|\le{n}} {\bf E}\{
G(i,i)G^{0}(p,p)\}\psi(\frac{i-p}{b}), \quad |i|\le{n} \eqno (3.19)
$$
where $G^{0}=\{G(i,j)-{\bf E}(G(i,j))\}_{ij}$ and $\epsilon_{ip}$ is
given by $(3.16)$. We are interested in the average value of the
normalized trace of the resolvent
$$
g_{n}(z)=N^{-1}\sum_{|i|\le{n}}G(i,i),
$$
where $G(i,i), \ |i|\le{n}$ are the diagonal entries of the resolvent. Taking the
normalized sum of $(3.15)$ and using the notation
$g_{n}(z)=g_{n,b}(z)$, we obtain
$$
{\bf E}\{g_{n,b}(z)\} = \xi + \frac{\xi v^{2}}{N}
\sum_{|i|,|p|\le{n}} {\bf E}G(i,i){\bf E}G(p,p)
\frac{\psi{(\frac{i-p}{b})}}{b}
$$
$$
+\{\frac{1}{N}\sum_{|i|\le{n}} \frac{1}{b} \phi_{1}{(i)}\}
+\{\frac{1}{N}\sum_{|i|\le{n}} \phi_{2}{(i)}\} -
\frac{\xi}{N}\sum_{|i|,|p|\le{n}} \epsilon_{ip}. \eqno (3.20)
$$
To prove $\lim{\bf E}\{g_{n,b}(z)\}=w(z)$ $(2.11)$, we need the
following statement concerning the pointwise convergence in average
of the diagonal entries $G(i,i)=G_{n,b}(i,i;z), \ |i|\le{n}$ of the
resolvent. Given a positive integer $L$, let
$$
B_{L}\equiv B_{L}(n,b)=\{ i\in \mathbf{ Z}\ : \ |i|\le{n-bL}\}.
\eqno (3.21)
$$

 \vskip0,5cm

{\bf Lemma 3.1} {\it Given  $\epsilon>0$, there exists a positive
integer $L=L(\epsilon)$ such that
$$
\sup_{i\in B_{L}}| {\bf E}\{G_{n,b}(i,i;z)\} - w(z)  |
\le{\epsilon}, \ z\in\Lambda_{\eta} \eqno (3.22)
$$
for sufficiently large $b$, $n$ satisfying $(2.3)$.

} \vskip0,5cm

 Lemma 3.1 will be proved in the next subsection. Let us assume that
$(3.22)$ is true. Then regarding definition of $B_{L}$,
we deduce from inequality $(3.7)$ that
$$
|{\bf E}\{g_{n,b}(z)\}-\frac{1}{N}\sum_{i\in B_{L}}{\bf
E}G(i,i)|=|\frac{1}{N}\sum_{i=-n}^{-n+bL-1}{\bf
E}G(i,i)+\frac{1}{N}\sum_{i=n-bL+1}^{n}{\bf E}G(i,i)|
$$

$$
\le{\frac{2bL-2}{N|Imz|}}, \quad z\in \Lambda_{\eta}. \eqno (3.23)
$$
Now $(2.11)$ follows from Lemma 3.1, inequality $(3.23)$ and the
limiting condition $(2.3)$. $\diamond$

\vskip0,5cm

\subsection{Proof of Lemma 3.1}
Let us have a step back and rewrite $(3.17)$ as
$$
{\bf E}G(i,i) = \xi+\xi v^{2}{\bf E}G(i,i){\bf E}U_{G}(i) +
\frac{1}{b} \phi_{1}{(i)} + \phi_{2}{(i)} - \xi\sum_{|p|\le{n}}
\epsilon_{ip}, \eqno (3.24)
$$
with
$$
U_{G}(i)=\frac{1}{b}\sum_{|p|\le{n}}G(p,p)\psi{(\frac{i-p}{b})},
$$
where $\phi_{1}$, $\phi_{2}$ and $\epsilon_{ip}$ given by relations
$(3.18)$, $(3.19)$ and $(3.16)$.
Let us denote the average ${\bf E}G(i,i)$ by $g(i)$ and rewrite $(3.24)$ in the following form:
$$
g(i) = \xi + \xi v^{2}g(i)U_{g}(i) + \frac{1}{b} \phi_{1}{(i)} +
\phi_{2}{(i)} - \xi\sum_{|p|\le{n}} \epsilon_{ip}.  \eqno (3.25)
$$
Let us consider the solution $\{r(i),\ |i|\le{n}\}$ of equation
$$
r(i) = \xi + \xi v^{2}r(i)U_{r}(i), \quad |i|\le{n}. \eqno (3.26)
$$
Given $z\in {\Lambda_{\eta}}$, one can prove that the system of
equation $(3.26)$ is uniquely solvable in the set of $N$-dimensional
vectors $\{\overrightarrow{r}\}$ such that
$$
|| \overrightarrow{r} || = sup_{|i|\le{n}} |r(i)|
\le{\frac{2}{|Im{(z)}|}}. \eqno (3.27)
$$
(see Lemma 6.1 of \cite{A}). Certainly, $r(i)$ and $g(i)$ $(3.25)$
depend on particular values of $z$, $n$ and $b$, so we shall use the
notations $r(i)=r_{n,b}(i;z)$ and $g(i)=g_{n,b}(i;z)$. The following
statements concern the differences:

$$
D_{n,b}(i;z)=g_{n,b}(i;z) - r_{n,b}(i;z) , \quad d_{n,b}(i;z)=
r_{n,b}(i;z) - w(z),
$$
where $w(z)$ verifies equation $(2.9)$. \vskip0,5cm

{\bf Lemma 3.2} {\it Given  $\epsilon>0$, there exists a positive
integer $L=L(\epsilon)$ such that for all sufficiently large $b$ and
$n$ satisfying $(2.3)$, inequality
$$
  \sup_{i\in B_{L}}|d_{n,b}(i;z)| \le{\epsilon}, \quad
  z\in\Lambda_{\eta},
\eqno (3.28)
$$
holds, with $B_{L}$ given by $(3.21)$.}

\vskip0,5cm

{\bf Lemma 3.3} {\it If $z\in\Lambda_{\eta}$ and $(2.3)$ holds, then
$$
\sup_{|p|\le{n}}|D_{n,b}(i,z)| = o(1), \quad n,b \longrightarrow
\infty,  \eqno (3.29)
$$
} \vskip0,5cm

 Lemma 3.1 follows from $(3.28)$ and $(3.29)$. Under the same conditions
of Lemma 3.1, one can find $L^{'}\geq{L}$ such that
$$
\sup_{i\in B_{L^{'}}}|\frac{\xi}{1-\xi v^{2}U_{g}(i)} -
w(z)|\le{2\epsilon}.
 \eqno (3.30)
$$
Relation $(3.30)$ follows from $(3.22)$ and $(3.27)$, a priori
estimate
$$
\sup_{|i|\le{n}}|g(i)|\le{\frac{1}{|Imz|}}, \eqno (3.31)
$$
and the observation that $L^{'}$ has to satisfy condition
$\psi(L-L^{'})\le{\epsilon}$.

\vskip0,5cm

{\it Proof of Lemma 3.2}. Let us consider function
$w_{i}(z)=w(z)$ satisfying $(2.9)$ that we rewrite in the following
form similar to $(3.26)$:
$$
w_{i}(z)=\xi + \xi
v^{2}w_{i}(z)\frac{1}{b}\sum_{t=1}^{n}b\delta_{it}w_{t}(z), \quad
|i|\le{n}.
$$
Subtracting this equality from $(3.26)$ , we see that $d(i)\equiv
d_{n,b}(i;z)$ is given by relation
$$
d(i)=\xi v^{2}d(i)U_{r}(i)+\xi v^{2}w(z)U_{d}(i)+\xi
v^{2}w^{2}(z)[P_{b}+T(i)], \eqno (3.32)
$$
where
$$
P_{b}=\frac{1}{b}\sum_{t\in
\Z}\psi(\frac{t}{b})-\int_{\mathbf{R}}\psi(s)ds \eqno (3.33)
$$
and
$$
T_{n,b}(i)\equiv T(i)=\frac{1}{b}\sum_{|t|\le{n}}\psi(\frac{t-i}{b})
- \frac{1}{b}\sum_{t\in \Z}\psi(\frac{t}{b}). \eqno (3.34)
$$
It is clear that $|P_{b}|=o(1)$ as $b\rightarrow \infty$. Indeed,
one can choose an even step-like function $\psi_{d}(t)$ , $t\in
\R$ such that
$$
\psi_{d}(t)=\sum_{k\in
\N}\psi(\frac{k}{b})I_{(\frac{k-1}{b},\frac{k}{b})}(t),
\quad t\geq{0}.
$$
Then $\psi_{d}(t)\le{\psi(t)}$, $\psi_{d}(t)\rightarrow \psi(t)$ as
$b\rightarrow \infty$ and the Beppo-L\'evy Theorem implies the
convergence of the corresponding integrals of $(3.33)$.

Using equality
$$
r(i)=\frac{\xi}{1-v^{2}\xi U_{r}(i)},
$$
we obtain relation
$$
d(i)=v^{2}wr(i)U_{d}(i)+v^{2}w^{2}r(i)[P_{b}+T_{n,b}(i)], \
|i|\le{n}, \eqno (3.35)
$$
where we write $w=w(z)$. This relation, together with estimates
$(3.27)$ and $|w(z)|\le{|Imz|^{-1}}$, implies that
$|d(i)|\le{2|Imz|^{-1}}<2$.

 Given $\epsilon >0$, let us find a number $Q$ that
$2\int_{Q}^{\infty}\psi(t)dt<\epsilon$. Setting
$\tau:=v^{2}\eta^{-2}<1/4$, we derive from $(3.35)$ such
$$
\sup_{i\in B_{L}}|d(i)|\le{\tau[\sup_{i\in B_{L-Q}}|d(i)|+\sup_{i\in
B_{L}}|T_{n,b}(i)|+P_{b}+2\epsilon+\frac{2}{b} ] },
$$
where we used the condition $(2.1)$ and the estimate
$$
\frac{1}{b}\sum_{s:|s-i|>Qb}|d(s)|\psi(\frac{s-i}{b})
\le{4\int_{Q}^{\infty}\psi(t)dt +\frac{2}{b}}
$$
that follows from the monotonicity of $\psi(t)$. Now, it is easy to
conclude that
$$
\sup_{i\in B_{L}}|d(i)|\le{\sum_{0\le{j}\le{L/Q}}\tau^{j}[\sup_{i\in
B_{L-jQ}}|T_{n,b}(i)|+P_{b}]+\tau^{L/Q}\sup_{|i|\le{n}}|d(i)|+4\epsilon+
\frac{4}{b} }. \eqno (3.36)
$$
Let us choose such $M$ that $\tau^{M}<\epsilon$, then
$$
\sum_{0\le{j}\le{L/Q}}\tau^{j}\sup_{i\in
B_{L-jQ}}|T_{n,b}(i)|\le{\sup_{i\in B_{L-MQ}}|T_{n,b}(i)|+2\epsilon
}.
$$
Finally, observe that
$$
\sup_{i\in
B_{L-MQ}}|T_{n,b}(i)|\le{\frac{2}{b}\sum_{t=n-(L-MQ)b}^{\infty}\psi(\frac{t}{b})
}\le{2\int_{n/b-L+MQ}^{\infty}\psi(s)ds+\frac{1}{b} }.
$$
This inequality shows that $(3.28)$ holds for sufficiently large $L$ and
$1<<b<<n$. $\diamond$

\vskip0,5cm

{\it Proof of Lemma 3.3}. Subtracting $(3.26)$ from $(3.25)$, we
obtain relation for $D(i)=D_{n,b}(i,z)$,
$$
D(i) = \xi v^{2}g(i)U_{D}(i) + \xi v^{2}D(i)U_{r}(i) + \frac{1}{b}
\phi_{1}(i) + \phi_{2}(i) - \xi\sum_{|p|\le{n}} \epsilon_{ip}.
$$
This relation can be written as
$$
\vec{D}=[ I - W^{(g,r)}]^{-1} \vec{\theta},
$$
where we denote by $W^{(g,r)}$ a linear operator acting on the
vectors $e=(e_{s})_{-n\le{s}\le{n}}$ as
$$
[W^{(g,r)}e](i) = v^{2}g(i)r(i)\sum_{|s|\le{n}}e(s)u(i,s)
$$
and
$$
\vec{\theta}^{(r)}_{n,b}(i) =\frac{r(i)}{\xi}  \{\frac{1}{b}
\phi_{1}{(i)} + \phi_{2}{(i)} - \xi\sum_{p=1}^{n} \epsilon_{ip} \}.
$$
It is easy to see that if $z\in{\Lambda_{\eta}}$, then the estimates
$(3.27)$ and $(3.31)$ imply
$$
||W^{(g,r)} ||\le{\frac{v^{2}}{|Imz|^{2}}}< \frac{1}{2}. \eqno
(3.37)
$$
Let us accept for the moment that
$$
\sup_{|i|\le{n}}|\frac{r(i)}{\xi}  \{\frac{1}{b} \phi_{1}{(i)} +
\phi_{2}{(i)} - \xi\sum_{|p|\le{n}} \epsilon_{ip} \}|=o(1), \quad
z\in{\Lambda_{\eta}} \eqno (3.38)
$$
in the limit $n,b\rightarrow\infty$. Then, Lemma 3.3 follows from
$(3.37)$ and estimate $(3.38)$. \vskip0,5cm

Now, let us prove $(3.38)$. Using inequality $(3.16)$, we obtain the
first estimate concerning $\epsilon_{ip}$,

$$
|\sum_{|p|\le{n}}\epsilon_{ip}|\le{\frac{8c\hat{\mu}_{3}}{|Imz|^{3}\sqrt{b}}
\sum_{|p|\le{n}}\frac{\psi(\frac{i-p}{b})}{b}
}=O(\frac{1}{\sqrt{b}}), \eqno (3.39)
$$
where we used condition $(2.1)$ and
$$
\sum_{|p|\le{n}}\frac{\psi(\frac{i-p}{b})}{b}\le{\sum_{t\in
\Z}\frac{\psi(t)}{b} }\le{2\int_{0}^{\infty}\psi(t)dt+
\frac{\psi(0)}{b}}.
$$

 The term $\phi_{1}$ in $(3.18)$ is estimated with the help of
the elementary inequality
$$
\sum_{|p|\le{n}}|G(i,p)|^{2} = ||G\vec{e}_{i}||^{2}
\le{\frac{1}{|Imz|^{2}}} \eqno (3.40)
$$
which follows from $(3.7)$ and the inequality
$||G^{2}(z)||\le{|Imz|^{-2}}$. Then, $(3.40)$ implies
$$
|\frac{1}{b} \phi_{1}{(i)}|\le{| \frac{\xi v^{2}}{b}{\bf
E}{(\sum_{|p|\le{n}}G(i,p)^{2})}|+|\frac{\xi v^{2}}{b}{\bf
E}\{G(i,p)^{2}\}\psi(0)| } \le{\frac{2v^{2}}{|Imz|^{3}b}}. \eqno
(3.41)
$$

Using relations $(2.1)$, $(3.7)$ and identity
$$
{\bf E}\{fg^{0}\}={\bf E}\{f^{0}g^{0}\},
$$
we obtain that
$$
|\phi_{2}(i)|=|\frac{\xi v^{2}}{b}\sum_{|p|\le{n}} {\bf E}\{
G(i,i)G^{0}(p,p) \}\psi{(\frac{i-p}{b})}| \le{
\frac{v^{2}}{|Imz|}{\bf E}|G^{0}(i,i)U^{0}_{G}(i)|}. \eqno (3.42)
$$

Thus, to prove $(3.29)$, it is sufficient to prove inequalities
$(3.39)$, $(3.41)$, $(3.42)$ and the estimate
$$
\sup_{|i|\le{n}}{\bf E}|G^{0}(i,i)U^{0}_{G}(i)|=o(1), \quad
z\in{\Lambda_{\eta}}. \eqno (3.43)
$$
We prove estimate $(3.43)$ in section 4. Assuming that this is done,
one can say that Lemma 3.1 and relation $(2.11)$ are proved.

\vskip0,5cm

\section{Proof of the main Theorems}
In this section, we study the variance ${\bf Var}\{g_{n,b}(z)\}$ and
complete the proof of Theorem 2.3. Finally, we prove Theorem 2.2 and
Theorem 2.1.

\vskip0,2cm

\subsection{The variance and the proof of Theorem 2.1}
 In section 3, we proved relation $(2.11)$ assuming that the asymptotic relation
$(3.43)$ holds. To complete the proof of Theorem 2.3, we
follow mainly two steps. The first step is to remark that the
estimate $(3.43)$ and $(2.12)$ are consequences of the fact that the
variance of the diagonal entries $G(i,i;z)$
$$
{\bf Var}\{G(i,i;z)\}={\bf E}|G^{0}(i,i;z)|^{2}
$$
vanishes as $n,b\rightarrow\infty$. The second step is to prove the
following two relations that concern the moments of diagonal
elements of $G$.

\vskip0,5cm

{\bf Lemma 4.1}. {\it If $z\in \Lambda_{\eta}$ $(2.10)$ and
$\hat{\mu}_{3}<\infty$ $(2.4)$, then the estimates
$$
\sup_{|i|\le{n}}{\bf E}|U^{0}_{G}(i)|^{2}=o(b^{-3/2}) \eqno (4.1)
$$
and
$$
\sup_{|i|\le{n}}{\bf E}|G^{0}(i,i;z)|^{2}=o(b^{-1/2}), \eqno (4.2)
$$
are true in the limit $n,b\rightarrow\infty$. }

\vskip0,5cm

It it easy to show that estimate $(3.43)$ follows from $(4.2)$.
Finally, $(2.12)$ follows from $(4.2)$ and inequalities
$$
{\bf Var}\{g_{n}(z)\} = {\bf E}|g_{n}{(z)} - {\bf
E}{(g_{n}{(z)})}|^{2}\le{ {\bf E}\{
\frac{1}{N^{2}}\sum_{|s|,|p|\le{n}}|G^{0}(p,p) G^{0}(s,s) |\} }
$$
$$
\le { \frac{1}{N^{2}}\sum_{|s|,|p|\le{n}} \sqrt{{\bf
E}{(|G^{0}(p,p)|^{2})}}\sqrt{{\bf E}{(|G^{0}(s,s)|^{2})}}} \le{
\sup_{|p|\le{n}}{\bf E}{(|G^{0}(p,p)|^{2})}}. \eqno (4.3)
$$
Theorem 2.3 is proved. $\diamond$ \vskip0,2cm

We close this subsection with the proof of Lemma 4.1.

 \vskip0,2cm

{\it Proof of Lemma 4.1.} Let us denote $G_{1}=(H-z_{1})^{-1}$,
$G_{2}=(H-z_{2})^{-1}$ with $z_{j}\in \Lambda_{\eta}, \ j=1,2$.
Regarding the average ${\bf E}\{G^{0}_{1}(i,i)G_{2}(p,p)\}$, we
apply to $G_{2}$ the resolvent identity and we use $(3.2)$, to get
$$
{\bf E}\{G^{0}_{1}(i,i)G_{2}(p,p)\}=-\xi_{2}\sum_{|s|\le{n}}{\bf
E}\{G^{0}_{1}(i,i)G_{2}(p,s)\}.
$$
For each pair $(s,p)$, expression $G^{0}_{1}(i,i)G_{2}(p,s)$
represents a smooth function of $H(s,p)$ and its derivative is
bounded because of $(3.6)$ and $(3.7)$. In particular
$|D^{2}_{sp}[G^{0}_{1}(i,i)G_{2}(p,s)]|\le{c_{1}(|Imz_{1}|^{-1}
+|Imz_{2}|^{-1})^{4}}$ where $c_{1}$ is a constant and
$D^{r}_{sp}=\partial^{r}/\partial H(s,p)^{r}$.

According to $(2.2)$, $(1.5)$ and condition $\hat{\mu}_{3}<\infty$
$(2.4)$, the third absolute moment of $H(s,p)$ is of the order
$b^{-3/2}$. Thus applying $(3.4)$ to ${\bf
E}\{G^{0}_{1}(i,i)G_{2}(p,s)\}$ with $q=1$, one obtains

$$
{\bf E}G^{0}_{1}(i,i)G_{2}(p,p) = \xi_{2} v^{2}{\bf E}\{
G^{0}_{1}(i,i)G_{2}(p,p)U_{G_{2}}(p)  \} + \xi_{2}
v^{2}\sum_{|s|\le{n}}{\bf E}\{G^{0}_{1}(i,i)G_{2}(p,s)^{2}u(s,p)  \}
$$
$$
+ 2\xi_{2} v^{2}\sum_{|s|\le{n}}{\bf
E}\{G_{1}(i,s)G_{1}(p,i)G_{2}(p,s)u(p,s)\}-\frac{\xi_{2}v^{2}}{b}{\bf
E}\{G^{0}_{1}(i,i)G_{2}(p,p)^{2}
$$
$$
+G_{1}(i,p)^{2}G_{2}(p,p)\}\psi(0)-
\xi_{2}\sum_{|s|\le{n}}\tau_{sp},
$$
where
$$
 |\tau_{sp}|\le{c \{\sup|D^{2}_{sp}[G^{0}_{1}(i,i)G_{2}(p,s)]|\} b^{-3/2}{\bf
E}{(|a(p,s)|^{3})}{\bf E}{(d(p,s)^{3})}} \eqno (4.4)
$$
with a constant c,
$$
 u(s,p)=\frac{\psi( \frac{s-p}{b}   )}{b}  \quad   and \quad
U_{G}(p)= \sum_{|s|\le{n}}G(s,s)u(s,p).
$$

 Let us introduce variables
$$
q_{2}(p) = \frac{\xi_{2}}{1-\xi_{2} v^{2}U_{{\bf E}(G_{2})}(p)}, \
|p|\le{n}. \eqno (4.5)
$$
 Using $(4.5)$ and identity
$$
{\bf E}\{fg\} = {\bf E}\{fg^{0}\} + {\bf E}\{f\}{\bf E}\{g\}, \eqno
(4.6)
$$
we get equality
$$
{\bf E}G_{1}^{0}(i,i)G_{2}(p,p) = q_{2}(p) v^{2}{\bf E}\{
G^{0}_{1}(i,i)G_{2}(p,p)U^{0}_{G_{2}}(p)  \}
$$
$$
+q_{2}(p)v^{2}\sum_{|s|\le{n}}{\bf
E}\{G^{0}_{1}(i,i)G_{2}(p,s)^{2}u(s,p)  \}
$$
$$
+ 2q_{2}(p)v^{2}\sum_{|s|\le{n}}{\bf
E}\{G_{1}(i,s)G_{1}(p,i)G_{2}(p,s)u(p,s) \}
$$
$$
-\frac{q_{2}(p)v^{2}}{b}{\bf
E}\{G^{0}_{1}(i,i)G_{2}(p,p)^{2}+G_{1}(i,p)^{2}G_{2}(p,p)\}\psi(0)
-q_{2}(p)\sum_{|s|\le{n}}\tau_{sp}. \eqno (4.7)
$$
Multiplying $(4.7)$ by $u(i,t)$ and summing over $i$, we get relation
$$
{\bf E}U^{0}_{G_{1}}(t)G_{2}(p,p) = q_{2}(p)v^{2}{\bf E}\{
U^{0}_{G_{1}}(t)G_{2}(p,p)U^{0}_{G_{2}}(p)  \}
$$
$$
+q_{2}(p)v^{2}\sum_{|s|\le{n}}{\bf
E}\{U^{0}_{G_{1}}(t)G_{2}(p,s)^{2}u(s,p)  \}
$$
$$
+ 2q_{2}(p)v^{2}\sum_{i,s}{\bf
E}\{G_{1}(i,s)G_{1}(p,i)G_{2}(p,s)u(p,s)u(i,t)   \}
$$
$$
-\frac{q_{2}(p)v^{2}}{b}{\bf
E}\{U^{0}_{G_{1}}(t)G_{2}(p,p)^{2}\}\psi(0)-\frac{q_{2}(p)v^{2}}{b}{\bf
E}\{G_{2}(p,p)\sum_{|i|\le{n}}G_{1}(i,p)^{2}u(i,t)\}\psi(0)
$$
$$
-\sum_{i,s}q_{2}(p)\tau_{sp}u(i,t). \eqno (4.8)
$$
Regarding $G_{1}(p,.)u(.,t)$ and $G_{2}(.,s)u(s,.)$ as
$n-$dimensional vectors, one gets
$$
|\sum_{i,s}{\bf E}\{G_{1}(i,s)G_{1}(p,i)G_{2}(p,s)\}u(p,s)u(i,t)|
$$
$$
\le{ ||G_{1}||\sqrt{ \sum_{|i|\le{n}}|G_{1}(p,i)u(i,t) |^{2}  }
\sqrt{\sum_{|s|\le{n}}|G_{2}(p,s)u(s,p) |^{2}} }. \eqno (4.9)
$$
 Inequalities $(3.7)$ and $(3.40)$ imply that the right-hand side of $(4.9)$
is bounded by $b^{-2}\eta^{-3}$. Regarding the last term of $(4.8)$
and using $(4.4)$ and $(3.7)$, we get
$$
|\sum_{|i|,|s|\le{n}} \tau_{sp}u(i,t)
|\le{\frac{cc_{1}}{16\eta^{4}\sqrt{b}} \hat{\mu}_{3}
[\sum_{|i|\le{n}}u(i,t)][\sum_{|s|\le{n}}u(p,s)]}=O(\frac{1}{\sqrt{b}}).
\eqno (4.10)
$$
The first term of the right-hand side of $(4.8)$ can be estimated by

$$
|q_{2}(p) v^{2}{\bf E}\{ U^{0}_{G_{1}}(t)G_{2}(p,p)U^{0}_{G_{2}}(p)
\}|\le{\frac{v^{2}}{|Imz_{2}|^{2}}\sqrt{\sup_{|i|\le{n}}{\bf
E}|U^{0}_{G_{1}}(i)|^{2}}\sqrt{\sup_{|i|\le{n}}{\bf
E}|U^{0}_{G_{2}}(i)|^{2}}}, \eqno (4.11)
$$
 where we used estimates $(3.7)$ and
$$
|q_{2}(p)|\le{\frac{1}{|Imz_{2}|}}. \eqno (4.12)
$$
Also, we use $(3.40)$ and $(4.12)$ to see that
$$
|q_{2}(p)v^{2}\sum_{|s|\le{n}}{\bf
E}\{U^{0}_{G_{1}}(t)G_{2}(p,s)^{2}u(s,p)
\}|\le{\frac{v^{2}}{b|Imz_{2}|^{3}}\sqrt{\sup_{|i|\le{n}}\{  {\bf
E}|U^{0}_{G_{1}}(i)|^{2}  \}} } \eqno (4.13)
$$
and
$$
|\frac{q_{2}(p)v^{2}}{b}{\bf
E}\{U^{0}_{G_{1}}(t)G_{2}(p,p)^{2}\}\psi(0)+\frac{q_{2}(p)v^{2}}{b}{\bf
E}\{G_{2}(p,p)\sum_{|i|\le{n}}G_{1}(i,p)^{2}u(i,t)\}\psi(0)|
$$
$$
\le{\frac{v^{2}}{b|Imz_{2}|^{3}}\sqrt{\sup_{|i|\le{n}}\{  {\bf
E}|U^{0}_{G_{1}}(i)|^{2}  \}}
+\frac{v^{2}}{b^{2}|Imz_{1}|^{2}|Imz_{2}|^{2}}}. \eqno (4.14)
$$

Now, multiplying  $(4.8)$ by $u(p,r)$ and summing  over  $p$, and
taking  $G_{1}=\bar{G}_{2}$ and  $i=r$,  we obtain a relation that
together  with  $(4.9)$, $(4.10)$, $(4.11)$, $(4.13)$ and $(4.14)$
implies the following estimate for the variable
$M_{12}=\sup_{|i|\le{n}}\{ {\bf E}|U^{0}_{G_{1}}(i)|^{2}  \}$:

$$
M_{12}\le{ \frac{v^{2}}{|Imz_{2}|^{2}}M_{12} +
\frac{A_{1}}{b}\sqrt{M_{12}} + \frac{A_{2}}{\sqrt{b}}},
$$
with $A_{1}$ and $A_{2}$ are constants. Regarding this inequality,
it is easy to show that $M_{12}=O(b^{-3/2})$. This proves $(4.1)$.

In order to prove $(4.2)$, we go back to relation $(4.7)$ with
$G_{1}=\bar{G}_{2}$,  $p=i$. Applying $(3.7)$ and $(4.12)$, we
obtain estimates
$$
\sup_{|i|\le{n}}|G^{0}_{1}(i,i)|^{2}\le{  \frac{v^{2}}{\eta^{3}}
\sqrt{{\bf E}(|U^{0}_{G_{1}}(i)|^{2})} +
\frac{v^{2}}{\eta^{2}b}\sum_{|s|\le{n}}|G_{1}(i,s)|^{2}    }
$$
$$
+ \frac{2v^{2}}{\eta^{2}b}\sum_{|s|\le{n}}|G_{1}(i,s)|^{2} +
\frac{3v^{2}}{\eta^{4}b}+ \frac{1}{\eta}\sum_{|s|\le{n}}|\tau_{si}|,
\eqno (4.15)
$$
with $\eta=2v+1$. Regarding the last term of $(4.15)$ and $(4.4)$,
we obtain inequality

$$
\sum_{|s|\le{n}}|\tau_{si}|\le{
\frac{c_{2}}{\sqrt{b}}\sum_{|s|\le{n}}u(s,i)
}=O(\frac{1}{\sqrt{b}}), \eqno (4.16)
$$
 where $c_{2}$ is a constant. Applying $(3.40)$, $(4.1)$ and $(4.16)$
to inequality $(4.15)$, we obtain $(4.2)$ and we are done.
$\diamond$

\subsection{Proof of Theorem 2.2 and Theorem 2.1}
{\it Proof of Theorem 2.2.}  It is easy to show that relation
$(2.11)$ together with $(2.12)$ implies convergence in
probability:
$$
P-\lim_{n\rightarrow\infty} g_{n}(z) = w(z),  \quad z\in
{\Lambda_{\eta}}.
$$
By definition $(2.6)$, we rewrite this convergence in the form
$$
\{\int_{\mathbf{R}}\frac{d\sigma_{n}(\lambda,H_{n})}{\lambda -
z}\}_{n\rightarrow\infty}\stackrel{P}{\longrightarrow}\{\int_{\mathbf{R}}
\frac{d\sigma(\lambda)}{\lambda- z}\}, \  z\in {\Lambda_{\eta}}.
\eqno (4.17)
$$
 Now we will prove that $(4.17)$ is true for all non real $z$.
Representation $(4.17)$ for $g_{n}(z)$ and $w(z)$ implies that
the function $z\mapsto{g_{n}(z) - w(z)}$ is analytic and uniformly
bounded on any compact set $\Gamma \subset{\Lambda_{\eta}}$. Thus,
given $\epsilon>0$, there exists a finite set $z_{1},...,z_{k}
\in{\Gamma}$ such that
$$
\max_{z\in{\Gamma}}|g_{n}(z) - w(z)| \le{ \epsilon +
\max_{j=1,..,k}|g_{n}(z_{j}) - w(z_{j})|}.
$$
 Let  $\{n_{1}\}$ a subsequence from $\{n\}$. We have  $g_{n_{1}}(z_{1})
\stackrel{P}{\longrightarrow}w(z_{1})$, then there exists a
subsequence $\{n_{1}^{'}\}\subseteq{\{n_{1}\}   }$ such that
$g_{n_{1}^{'}}(z_{1})\stackrel{a.s.}{\longrightarrow}w(z_{1})$.
Therefore, the sequence $\{n_{1}^{'}\}$ contains a subsequence
$\{n_{2}^{'}\}$ where
$g_{n_{2}^{'}}(z_{2})\stackrel{a.s.}{\longrightarrow}w(z_{2})$ and
$g_{n_{2}^{'}}(z_{1})\stackrel{a.s.}{\longrightarrow}w(z_{1})$. Now we see that there exist a subsequence
$n^{'}=n_{k}^{'}\subseteq{...}\subseteq{n_{1}^{'}}$ common to all
$z_{1},...,z_{k}$ such that
$$ \max_{j=1,..,k}|g_{n^{'}}(z_{j}) -
w(z_{j})|\stackrel{a.s.}{\longrightarrow}0.
$$
As a result, we get consequence
 $$
\lim_{n^{'}\rightarrow+\infty}\{ \max_{z\in{\Gamma}}|g_{n^{'}}(z) -
w(z)|\} =0  \quad   a.s..
$$
Now we obtain $g_{n^{'}}(z)\stackrel{a.s.}{\longrightarrow}w(z)$
with $\{n^{'}\}$ is common for all $z\in{\Gamma}$. Since
$g_{n^{'}}(z)$ and $w(z)$ are analytic functions in $ \mathbb{C}
\setminus \mathbb{R} $, and
$g_{n^{'}}(z)\stackrel{a.s.}{\longrightarrow}w(z)$ for
$z\in{\Gamma}$, then this convergence is true for all $ z\in
{\mathbb{C} \setminus \mathbb{R} }$. Thus, we rewrite $(4.16)$ as
$$
\{\int_{\mathbf{R}}\frac{d\sigma_{n^{'},b^{'}}(\lambda,H_{n^{'}})}{\lambda
-
z}\}_{n^{'},b^{'}\rightarrow+\infty}\stackrel{a.s.}{\longrightarrow}\{\int_
{\mathbf{R}}\frac{d\sigma(\lambda)}{\lambda- z}\}, \quad  z\in
{\Lambda_{0}}. \eqno (4.18)
$$

The sequence $\{n_{1}\}$ contains a subsequence $\{n^{'}\}$ such
that $(4.18)$ holds, which implies convergence in probability
$(2.9)$ (see Lemma in paper \cite{C}).  $\diamond$

\vskip0,5cm

{\it Proof of Theorem 2.1.} It is well known that the set of linear
combinations of $\{\frac{1} {\lambda - z_{j}} ; \
z_{j}\in\Lambda_{0}\}$ is dense in ${C_{0}(\mathbf{R})}$ (set of
counting functions), then $(4.18)$ implies that
$$
\{\int_{\mathbf{R}}\varphi(\lambda)d\sigma_{n^{'}}(\lambda,H_{n^{'}})\}_{n^{'}\rightarrow
\infty}\stackrel{a.s.}{\longrightarrow}\{\int_{\mathbf{R}}\varphi(\lambda)
d\sigma(\lambda)\} \eqno (4.19)
$$
for all $\varphi\in{C_{0}(\mathbf{R})}$.

The sequence $\{n_{1}\}$ contains a subsequence as $\{n^{'}\}$ such
that the convergence $(4.19)$ is true. We use the same arguments
as in the proof of Theorem 2.2, the convergence in probability
$(2.5)$ gives and finishes the proof of Theorem 2.1. $\diamond$
\vskip0,5cm

\section{Estimates of the Variance }
In this section, we obtain more estimates of the variance under more
restrictive conditions on the moments of $a(i,p)$ and the function
$\psi(t)$ given by $(2.1)$. Finally, we provide a proof of the
cumulant expansions formula.\vskip0,5cm

\subsection{First Estimate}
{\bf Theorem 5.1}. {\it Assume that ${\bf E}\{a^{2l+1}(i,p)\}= 0 $
with $l=0,1,2$ and
$$
{\bf E}\{a(i,p)^{2}\}=v^{2},\ {\bf E}\{a(i,p)^{4}\}=3v^{4}, \ {\bf
E}\{a(i,p)^{6}\}=15v^{6}, \quad |i|,|p|\le{n},
$$
$\hat{\mu}_{7}<\infty$ $(2.4)$, then the estimate
$$
{\bf Var}\{g_{n,b}(z)\}=O(\frac{1}{b})  \eqno (5.1)
$$
holds for large enough $n$ and $b$ and for all $z\in\Lambda_{\eta}$.
} \vskip0,5cm

{\it Proof of Theorem 5.1.} We prove Theorem 5.1 by using the
following estimates of the moments of the diagonal elements of the
resolvent $G$
 \vskip0,5cm

{\bf Lemma 5.1}. {\it Under conditions of Theorem 5.1, the estimates
$$
\sup_{|i|\le{n}}{\bf E}(|U^{0}_{G}(i;z)|^{2})=O(\frac{1}{b^{2}})
\eqno (5.2)
$$
and
$$
\sup_{|i|\le{n}}{\bf E}(|G(i,i;z)|^{2}) = O(\frac{1}{b}) \eqno (5.3)
$$
hold in the limit  $n,b\longrightarrow \infty$ and for all
$z\in\Lambda_{\eta}$. }

\vskip0,5cm

Now it is easy to show that Theorem 5.1 follows from inequality
$(4.3)$ and estimate $(5.3)$. \vskip0,5cm

{\it Proof of Lemma 5.1}. We start with $(5.2)$. Let us denote
$G_{1}=(H-z_{1})^{-1}$, $G_{2}=(H-z_{2})^{-1}$ with $z_{j}\in
\Lambda_{\eta}, \ j=1,2$. Consider the average ${\bf
E}U^{0}_{G_{1}}(t)G_{2}(p,p)$ and use $(3.2)$ to obtain
$$
{\bf E}U^{0}_{G_{1}}(t)G_{2}(p,p)=-\xi\sum_{s=1}^{n}{\bf
E}\{U^{0}_{G_{1}}(t)G_{2}(p,s)\}.
$$
For each pair $(s,p)$ $U^{0}_{G_{1}}(t)G_{2}(p,s)$ is a smooth
function of $H(s,p)$ and its derivative is bounded because of
$(3.7)$. In particular
$|D^{6}_{sp}[U^{0}_{G_{1}}(t)G_{2}(p,s)]|\le{c_{3}(|Imz_{1}|^{-1}
+|Imz_{2}|^{-1})^{8}}$ for some constant $c_{3}$ and
$D^{r}_{sp}=\partial^{r}/\partial H^{r}(s,p)$.

According to $(2.2)$, $(1.5)$ and condition $\hat{\mu}_{7}<\infty$
$(2.4)$, the 7-th absolute moment of $H(s,p)$ is of order
$b^{-7/2}$. Thus, applying $(3.4)$ to ${\bf
E}U^{0}_{G_{1}}(t)G_{2}(p,s)$ with $q=5$, one obtains

$$
{\bf E}\{U^{0}_{G_{1}}(t)G_{2}(p,p)\} =
-\xi_{2}\sum_{|s|\le{n}}K_{2}{\bf
E}\{D^{1}_{sp}[U^{0}_{G_{1}}(t)G_{2}(p,s)]\}
$$

$$
-\xi_{2}\sum_{|s|\le{n}}\frac{K_{4}}{6}{\bf E}\{
D^{3}_{sp}[U^{0}_{G_{1}}(t)G_{2}(p,s)]\} -\xi_{2}
\sum_{|s|\le{n}}\frac{K_{6}}{120}{\bf
E}\{D^{5}_{sp}[U^{0}_{G_{1}}(t)G_{2}(p,s)]\}
$$
$$
-\xi_{2}\sum_{|s|\le{n}}\epsilon^{(1)}_{sp}
$$
where
$$
|\epsilon^{(1)}_{sp}|\le{c\sup|D^{6}_{sp}[U^{0}_{G_{1}}(t)G_{2}(p,s)]|b^{-7/2}{\bf
E}(|a_{sp}|^{7}){\bf E}(d_{sp}^{7})}. \eqno (5.4)
$$
Regarding the terms involving the first derivative and using
$(3.6)$, the definition of $K_{2}$ $(3.13)$, we get
$$
{\bf E}\{U^{0}_{G_{1}}(t)G_{2}(p,p)\} =  \xi_{2}v^{2}{\bf
E}\{U^{0}_{G_{1}}(t)G_{2}(p,p)U_{G_{2}}(p)  \}
$$
$$
+\xi_{2}v^{2}\sum_{|s|\le{n}}{\bf
E}\{U^{0}_{G_{1}}(t)G_{2}(p,s)^{2}\}u(p,s)
$$
$$
+ 2\xi_{2}v^{2}\sum_{s,i}{\bf
E}\{G_{1}(i,s)G_{1}(p,i)G_{2}(p,s)\}u(p,s)u(i,t) + R(p),
$$
where
 $$
 R(p)=-\frac{\xi_{2}v^{2}}{b}{\bf
E}\{U^{0}_{G_{1}}(t)G_{2}(p,p)^{2}+
\sum_{|i|\le{n}}G_{1}(i,p)^{2}G_{2}(p,p)u(i,t)\}\psi(0)
$$
$$
-\xi_{2}\sum_{|s|\le{n}}\frac{K_{4}}{6}{\bf E}
 \{D^{3}_{sp}[U^{0}_{G_{1}}(t)G_{2}(p,s)]\}
-\xi_{2}\sum_{|s|\le{n}}\frac{K_{6}}{120}{\bf E}
 \{D^{5}_{sp}[U^{0}_{G_{1}}(t)G_{2}(p,s)]\}
$$
$$
-\xi_{2}\sum_{|s|\le{n}}\epsilon^{(1)}_{sp}.
$$
This term does not contribute to the estimate $(5.1)$.

   Applying $(4.6)$ to the first term and using the definition
of $q_{2}(p)$ $(4.5)$, one gets

$$
{\bf E}\{U^{0}_{G_{1}}(t)G_{2}(p,p)\}=q_{2}(p)v^{2}{\bf
E}\{U^{0}_{G_{1}} (t)G_{2}(p,p)U^{o}_{G_{2}}(p)  \}
$$
$$
+q_{2}(p)v^{2}\sum_{|s|\le{n}}{\bf
E}\{U^{0}_{G_{1}}(t)G_{2}(p,s)^{2}u(s,p)  \}
$$
$$
+ 2q_{2}(p)v^{2}\sum_{s,i}{\bf
E}\{G_{1}(i,s)G_{1}(p,i)G_{2}(p,s)u(p,s)u(i,t)   \}+
\frac{q_{2}(p)}{\xi_{2}}R(p). \eqno (5.5)
$$

Inequality $(4.9)$ implies that the third term of $(5.5)$ is of
order $O(b^{-2})$ as $n,b\longrightarrow \infty$.

Let us estimate each part of $(5.5)$ with the help of $U^{0}_{G}$.
Using $(3.40)$ and inequality $(4.12)$, we can estimate the first
term of the right-hand side of $(5.5)$ by

$$
|q_{2}(p)v^{2}{\bf E}\{ U^{0}_{G_{1}}(t)G_{2}(p,p)U^{0}_{G_{2}}(p)
\}|\le{\frac{v^{2}}{|Imz_{2}|^{2}} \sqrt{\sup_{|i|\le{n}}  {\bf
E}|U^{0}_{G_{1}}(i)|^{2} \sup_{|i|\le{n}}{\bf
E}|U^{0}_{G_{2}}(i)|^{2}} }. \eqno (5.6)
$$
Also we can write that
$$
|q_{2}(p)v^{2}\sum_{|s|\le{n}}{\bf
E}\{U^{0}_{G_{1}}(t)G_{2}(p,s)^{2}u(s,p)
\}|\le{\frac{v^{2}}{b|Imz_{2}|^{3}}\sqrt{\sup_{|i|\le{n}}\{  {\bf
E}|U^{0}_{G_{1}}(i)|^{2}\}} }. \eqno (5.7)
$$
If one assumes for a while that
$$
\sup_{|p|\le{n}}|R(p)|=O\{\frac{1}{b^{2}} + \frac{1}{b}{\bf
E}(|U^{0}_{G_{1}}(t)|)\}, \eqno (5.8)
$$
holds, then $(5.2)$ will follow. Indeed, multiplying $(5.5)$ by
$u(p,r)$ and summing over $p$, and taking $G_{1} = \bar{G}_{2}$ and
$i=r$, we obtain a relation that together with $(4.9)$, $(5.6)$,
$(5.7)$ and $(5.8)$ imply the following estimate for variable
$M_{12} = \sup_{|t|\le{n}}|U^{0}_{G_{1}}(t)|^{2}$:
$$
M_{12}\le{ \frac{v^{2}}{\eta^{2}}M_{12}+\frac{A_{3}}{b}\sqrt{M_{12}}
+ \frac{A_{4}}{b^{2}}},
$$
with $A_{3}$ and $A_{4}$ are constants. Regarding this inequality,
it is easy to show that $M_{12}=O(b^{-2})$. This implies $(5.2)$.
\vskip0,5cm

Let us prove $(5.8)$. Using inequality $(5.4)$, we get
$$
|\sum_{|s|\le{n}}\epsilon^{(1)}_{sp}|\le{
\frac{cc_{3}}{(2\eta)^{8}b^{5/2}}\sum_{|s|\le{n}}u(s,p)}=O(\frac{1}{b^{5/2}}), \eqno (5.9)
$$
for some constants $c$ and $c_{3}$. By the definition of the
cumulants (see subsection 3.1), we obtain
 $$
 K_{4}(Y) = {\bf E}(Y^{4}) - 3{\bf E}(Y^{2})^{2}=
 \frac{v^{4}}{b^{2}}  \Delta_{sp},
\eqno (5.10)
$$
where $Y=H(s,p)$, $ \Delta_{ip} =
3\psi{(\frac{s-p}{b})}\{1-\psi{(\frac{s-p}{b})}\}$ and
$$
 K_{6}(Y) = {\bf E}(Y^{6})-15{\bf E}(Y^{4}){\bf E}(Y^{2})
 - 10{\bf E}(Y^{3})^{2} + 30{\bf E}(Y^{2})^{3}
= \frac{v^{6}}{b^{3}}\Theta_{sp} \eqno (5.11)
$$
with $ \Theta_{ip} = 15\psi(\frac{s-p}{b})\{1
-\frac{45}{15}\psi(\frac{s-p}{b}) +
\frac{30}{15}\psi(\frac{s-p}{b})^{2}\}$.

 Using $(3.6)$ and $(3.7)$, we obtain that
$\sup_{t,s,p}|D^{r}_{sp}[U^{0}_{G_{1}}(t)G_{2}(p,s)]|=O(1)$ with
$r=1,2,..$. Which gives
$$
|\sum_{|s|\le{n}}\frac{K_{6}}{120}{\bf
E}\{D^{5}_{sp}[U^{0}_{G_{1}}(t)G_{2}(p,s)]\}| \le{
\frac{M_{1}}{b^{2}}\sum_{|s|\le{n}}\frac{\Theta_{sp}}{b}}
$$
$$
\le{\frac{15M_{1}}{b^{2}}\sum_{|s|\le{n}}\frac{\psi(\frac{s-p}{b})}{b}
}=O(\frac{1}{b^{2}}) \eqno (5.12)
$$
where $M_{1}$ is a constant.

One can estimate the terms with the third derivative by using
similar estimates in relations $(4.9)$ and $(5.6)$:
$$
-\xi_{2}\sum_{|s|\le{n}}\frac{K_{4}}{6}{\bf
E}\{D^{3}_{sp}[U^{0}_{G_{1}}(t)G_{2}(p,s)]\} = O\{\frac{1}{b_{2}} +
\frac{1}{b}{\bf E}(|U^{0}_{G_{1}}(t)|)\}. \eqno (5.13)
$$
Indeed,  ${\bf E}\{D^{3}_{sp}[U^{0}_{G_{1}}(t)G_{2}(p,s)]\} $
contains $44$ terms of the form
$$\sum_{|i|\le{n}} \  {\bf E}\{G_{\gamma_{1}}(\alpha_{1},\beta_{1})G_{\gamma_{2}}
(\alpha_{2},\beta_{2})G_{\gamma_{3}}(\alpha_{3},\beta_{3})G_{\gamma_{4}}(\alpha_{4},\beta_{4})G_{\gamma_{5}}(\alpha_{5},\beta_{5})\}
u(i,t) ,
$$
and $3$ terms of the form
$$
\sum_{|i|\le{n}} \  {\bf
E}\{G^{0}_{1}(i,i)G_{\gamma_{2}}(\alpha_{2},\beta_{2})G_{\gamma_{3}}
(\alpha_{3},\beta_{3})G_{\gamma_{4}}(\alpha_{4},\beta_{4})G_{\gamma_{5}}
(\alpha_{5},\beta_{5})\} u(i,t),
$$
where  $\gamma_{j}=1$ or $2$ and $\alpha_{j}\in\{i,s,p\}$ and
$\beta_{j}\in\{i,s,p\}$. These terms can be gathered into three
groups. In each group the terms are estimated by the same values
with the help of the same computations.

 We give estimates for the typical cases.
Using  $(2.1)$ and $(3.40)$, we get for the first terms of the first
group:
$$
|\sum_{i,s}K_{4}(H(s,p)) {\bf
E}\{G_{1}(i,p)^{2}G_{1}(s,s)G_{1}(p,s)G_{2}(p,s) \}u(i,t) |
$$
$$
\le{\frac{3^{4}}{b^{2}\eta^{3}}\sum_{|i|\le{n}}|G_{1}(i,p)^{2}|\sum_{|s|\le{n}}
u(s,p)}\le{ \frac{3v^{4}}{b^{2}\eta^{5}} }. \eqno (5.14)
$$
For the terms of the second group, we obtain estimates

$ |\sum_{i,s}K_{4}{\bf
E}\{G_{1}(i,s)G_{1}(i,p)G_{1}(p,p)G_{1}(s,s)G_{2}(p,s)\}u(i,t)| $
$$
\le{\frac{3v^{4}}{b\eta^{3}}\sum_{s,i}{\bf
E}\{|G_{1}(p,i)G_{2}(p,s)|\}u(p,s)u(i,t)}
$$
$$
\le{ \frac{3v^{4}}{b|Imz|^{3}}\sqrt{
\sum_{|i|\le{n}}|G_{1}(p,i)u(i,t) |^{2} }
\sqrt{\sum_{|s|\le{n}}|G_{2}(p,s)u(s,p) |^{2}} }=O(\frac{1}{b^{3}}).
\eqno (5.15)
$$
Finally, for the terms of the third group, we get inequalities
$$
|\sum_{i,s}K_{4}{\bf
E}\{G_{1}(i,i)^{0}G_{2}(p,s)^{2}G_{2}(p,p)G_{2}(s,s)\}u(i,t)|
$$
$$
\le{ \frac{3v^{4}}{b\eta^{4}}{\bf
E}(|U^{0}_{G_{1}}(t)|)\{\sum_{|s|\le{n}}u(s,p)\}     }
=O(\frac{1}{b}{\bf E}(|U^{0}_{G_{1}}(t)|)). \eqno (5.16)
$$

Gathering all the estimates of $47$ terms, we obtain $(5.13)$.
Finally, notice that the term
$$
|\frac{\xi_{2}v^{2}}{b}{\bf E}\{U^{0}_{G_{1}}(t)G_{2}(p,p)^{2}+
\sum_{|i|\le{n}}G_{1}(i,p)^{2}G_{2}(p,p)u(i,t)\}\psi(0)|
$$
is estimated using $(4.14)$. Then $(5.8)$ follows from $(4.14)$,
$(5.9)$, $(5.12)$ and $(5.13)$. $\diamond$ \vskip1cm

   Let us prove $(5.3)$. Using $(3.2)$, one has
$$
{\bf E}G^{0}_{1}(i,i)G_{2}(p,p)=-\xi_{2}\sum_{|s|\le{n}}{\bf
E}\{G^{0}_{1}(i,i)G_{2}(p,s)\}.
$$
For each pair $(s,p)$ $G^{0}_{1}(i,i)G_{2}(p,s)$ is a smooth
function of $H(s,p)$ and its derivative are bounded because of
$(3.6)$ and $(3.7)$. In particular
$|D^{4}_{sp}[G^{0}_{1}(i,i)G_{2}(p,s)]|\le{c_{4}(|Imz_{1}|^{-1}
+|Imz_{2}|^{-1})^{6}}$ with a constant $c_{4}$.

According to $(2.2)$, $(1.5)$ and condition $\hat{\mu}_{5}<\infty$
$(2.4)$, the 5-th absolute moment of $H(s,p)$ is of order
$b^{-3/2}$. Thus applying $(3.4)$ to ${\bf
E}\{G^{0}_{1}(i,i)G_{2}(p,s)\}$ with $q=3$, one obtains

$$
{\bf E}G_{1}^{0}(i,i)G_{2}(p,p) = \xi_{2}v^{2}{\bf E}\{
G^{0}_{1}(i,i)G_{2}(p,p)U_{G_{2}}(p)  \}
$$
$$
+ \xi_{2}v^{2}\sum_{|s|\le{n}}{\bf
E}\{G^{0}_{1}(i,i)G_{2}(p,s)^{2}u(s,p)  \} +
2\xi_{2}v^{2}\sum_{|s|\le{n}}{\bf
E}\{G_{1}(i,s)G_{1}(p,i)G_{2}(p,s)u(p,s)   \}
$$
$$
-\xi_{2}\sum_{|s|\le{n}}\frac{K_{4}}{6}{\bf
E}\{D^{3}_{sp}[G^{0}_{1}(i,i)G_{2}(p,s)]\}
-\xi_{2}\sum_{|s|\le{n}}\epsilon^{(2)}_{sp}
$$
with
$$
 |\epsilon^{(2)}_{sp}|\le{c \{\sup|D^{4}_{sp}[G^{0}_{1}(i,i)G_{2}(p,s)]|\}
b^{-5/2}{\bf E}{(|a_{ps}|^{5})}{\bf E}{(d_{ps}^{5})}}.
 $$
Using the identity $(4.6)$ and the definition of $q_{2}$  $(4.5)$,
we write
$$
{\bf E}G_{1}^{0}(i,i)G_{2}(p,p) = q_{2}(p)v^{2}{\bf E}\{
G^{0}_{1}(i,i)G_{2}(p,p)U^{0}_{G_{2}}(p)  \}
$$
$$
+q_{2}(p)v^{2}\sum_{|s|\le{n}}{\bf
E}\{G^{0}_{1}(i,i)G_{2}(p,s)^{2}u(s,p)  \}
$$
$$
+ 2q_{2}(p)v^{2}\sum_{|s|\le{n}}{\bf
E}\{G_{1}(i,s)G_{1}(p,i)G_{2}(p,s)u(p,s) \}
$$
$$
-\frac{q_{2}(p)v^{2}}{b}{\bf
E}\{G^{0}_{1}(i,i)G_{2}(p,p)^{2}+G_{1}(i,p)^{2}G_{2}(p,p)\}\psi(0)
$$
$$
-q_{2}(p)\sum_{|s|\le{n}}\frac{K_{4}}{6}{\bf
E}\{D^{3}_{sp}[G^{0}_{1}(i,i)G_{2}(p,s)]\} -
q_{2}(p)\sum_{|s|\le{n}}\epsilon^{(2)}_{sp}. \eqno (5.17)
$$
It is easy to show that
$$
|\sum_{|s|\le{n}}\epsilon^{(2)}_{sp} |\le{\frac{M_{2}}{b\sqrt{b}}
\sum_{s}u(s,p)    }=O(\frac{1}{b\sqrt{b}}), \eqno (5.18)
$$
with $M_{2}$ some constant. Inequality $(3.7)$ implies that
$D^{r}_{sp}[G_{1}^{0}(i,i)G_{2}(p,s)]=O(1)$ with $r=1,2,..$, whence
one obtains inequality

$$
 |\sum_{|s|\le{n}}K_{4}{\bf
E}\{D^{3}_{sp}[G_{1}^{0}(i,i)G_{2}(p,s)]\} |\le{
\frac{3v^{4}M_{3}}{b}\sum_{|s|\le{n}}u(s,p) }, \eqno (5.19)
$$
for some constant $M_{3}$. Finally, it is easy to show that $(3.7)$
implies estimate
$$
|\frac{q_{2}(p)v^{2}}{b}{\bf
E}\{G^{0}_{1}(i,i)G_{2}(p,p)^{2}+G_{1}(i,p)^{2}G_{2}(p,p)\}\psi(0)|
=O(\frac{1}{b}). \eqno (5.20)
$$

 Now, $(5.3)$ follows from relation $(5.17)$, inequality $(3.40)$ and
estimates $(5.2)$, $(5.18)$, $(5.19)$ and $(5.20)$. The estimate
$(5.3)$ is proved, then Lemma 5.1 follows and finishes the proof of
Theorem 5.1. $\diamond$

\vskip0,5cm

\subsection{Second Estimates for the Variance}

 In this subsection, we derive estimates of the variance under more strong
conditions on $H_{n}$ and $\psi$.

\vskip0,5cm

{\bf Theorem 5.2}. {\it Suppose that ${\bf E}\{a^{2l+1}(i,p)\}= 0 $
with $l=0,1,2$,
$$
{\bf E}\{a(i,p)^{2}\}=v^{2},\ {\bf E}\{a(i,p)^{4}\}=3v^{4}, \ {\bf
E}\{a(i,p)^{6}\}=15v^{6}, \quad |i|,|p|\le{n}
$$
and
$$
\hat{\mu}_{10}=\sup_{i,p}{\bf E}|a(i,p)|^{10}<\infty \quad and \quad
\int_{\mathbf{R}}\sqrt{\psi(t)}dt\le{\infty},
$$
then the estimates
$$
\sup_{|i|\le{n}}|{\bf E}g_{n}^{0}(z)G^{0}(i,i)| = O(\frac{1}{nb} +
\frac{1}{b}\sqrt{{\bf Var}g_{n}(z)})  \eqno (5.21)
$$
and
$$
{\bf Var}\{g_{n}(z)\}=O(\frac{1}{b^{2}}) \eqno (5.22)
$$
hold for large enough $n$ and $b$ and for all $z\in\Lambda_{\eta}$. } \vskip0,5cm

{\it Proof of Theorem 5.2}. We start with $(5.21)$. The proof
follows the lines of the proof of $(5.3)$. Regarding $(5.17)$ and
summing it over $i$ and using the notation $g_{n}(z_{1})=g_{1}$, we
obtain that
$$
N{\bf E}\{g^{0}_{1}G_{2}(p,p)\} = Nq_{2}(p)v^{2}{\bf E}\{
g^{0}_{1}G_{2}(p,p)U^{0}_{G_{2}}(p)  \}
$$
$$
+Nq_{2}(p)v^{2}\sum_{|s|\le{n}}{\bf
E}\{g^{0}_{1}G_{2}(p,s)^{2}u(s,p) \} + 2q_{2}(p)v^{2}\sum_{s,i}{\bf
E}\{G_{1}(i,s)G_{1}(p,i)G_{2}(p,s)u(p,s)\}
$$
$$
-\frac{q_{2}(p)v^{2}}{b}{\bf
E}\{Ng_{1}^{0}G_{2}(p,p)^{2}+\sum_{|i|\le{n}}G_{1}(i,p)^{2}G_{2}(p,p)\}\psi(0)
$$
$$
-Nq_{2}(p)\sum_{|s|\le{n}}\frac{K_{4}}{6}{\bf
E}\{D^{3}_{sp}[g_{1}^{0}G_{2}(p,s) ]\}
-q_{2}(p)\sum_{s,i}\epsilon^{(2)}_{sp}. \eqno (5.23)
$$
Using $(5.2)$, $(3.7)$, $(3.40)$ and $(4.12)$, it is easy to show
that
$$
Nq_{2}(p)v^{2}{\bf E}\{ g^{0}_{1}G_{2}(p,p)U^{0}_{G_{2}}(p)
\}=O(\frac{N}{b}\sqrt{{\bf Var}(g_{1})}), \eqno (5.24)
$$
$$
Nq_{2}(p)v^{2}\sum_{|s|\le{n}}{\bf E}\{g^{0}_{1}G_{2}(p,s)^{2}u(s,p)
\}=O(\frac{N}{b}\sqrt{{\bf Var}(g_{1})}), \eqno (5.25)
$$

$$
|2q_{2}(p)v^{2}\sum_{s,i}{\bf
E}\{G_{1}(i,s)G_{1}(p,i)G_{2}(p,s)u(p,s)   \} |
$$
$$
\le{\frac{2v}{|Imz_{1}||Imz_{2}|b}\sqrt{ \sum_{|i|\le{n}}|G_{1}(p,i)
|^{2}  } \sqrt{\sum_{|s|\le{n}}|G_{2}(p,s) |^{2}}}= O(\frac{1}{b})
 \eqno (5.26)
$$
and
$$
\frac{q_{2}(p)v^{2}}{b}{\bf
E}\{Ng_{1}^{0}G_{2}(p,p)^{2}+\sum_{|i|\le{n}}G_{1}(i,p)^{2}G_{2}(p,p)\}\psi(0)
$$
$$
=O(\frac{N}{b}\sqrt{{\bf Var}(g_{1})}+\frac{1}{b}). \eqno (5.27)
$$
Now, regarding $D^{r}_{sp}[g_{1}^{0}G_{2}(p,s)]$ with $r=1,2,..$ and
$(5.31)$ (see Lemma 5.2 at the end of this subsection), it is easy
to see that
$$
\sup_{s,p}|D^{r}_{sp}[g_{1}^{0}G_{2}(p,s)]| =
O(\frac{1}{N}+|g^{0}_{1}|), \quad r=1,2,...   \eqno (5.28)
$$
Then we obtain
$$
|Nq_{2}(p)\sum_{|s|\le{n}}K_{4}{\bf
E}\{D^{3}_{sp}[g_{1}^{0}G_{2}(p,s)]\}|
\le{\frac{N}{\eta}\sum_{|s|\le{n}}\frac{\Delta_{sp}}{b^{2}}   |{\bf
E}\{D^{3}_{sp}[g_{1}^{0}G_{2}(p,s)]\}|}
$$
$$
 = O(\frac{N}{b}[\frac{1}{N}+\sqrt{{\bf Var}(g_{1})}]).
\eqno (5.29)
 $$

Now we give more details for the remainder $\epsilon^{(2)}$

$$
|\epsilon^{(2)}_{sp}|\le{c\{K_{4}{\bf E}|H(s,p)
D^{4}_{sp}[G^{0}_{1}(i,i)G_{2}(p,s)]^{(1)} | +K_{2}{\bf
E}|H^{3}(s,p)D^{4}_{sp}[G^{0}_{1}(i,i)G_{2}(p,s)]^{(2)} | }
$$
$$
+ {\bf E}|H^{5}(s,p)D^{4}_{sp}[G^{0}_{1}(i,i)G_{2}(p,s)]^{(3)}|\},
\eqno (5.30a)
$$
where for each pair $(p,s)$

$$
[G^{0}_{1}(i,i)G_{2}(p,s)]^{(\nu)}=\{G^{(\nu)}\}^{0}_{sp}(i,i;z_{1})G^{(\nu)}_{sp}
(p,s;z_{2}), \quad \nu=1,2,3
$$
and $G^{(\nu)}_{sp}(z_{j})=(H^{(\nu)}_{sp}-z_{j})^{-1}$, $j=1,2$ is
the resolvent of the real symmetric matrix $H^{(\nu)}_{sp}$:
 $$
  H^{(\nu)}_{sp}(r,i)= \left\{
\begin{array}{lll}
H(r,i) & \textrm{if} & (r,i)\neq(s,p) \\
H^{(\nu)}(s,p) &  \textrm{if} & (r,i)=(s,p)
\end{array}\right.
$$
with $|H^{(\nu)}(s,p)|\le{|H(s,p)|}$, $\nu=1,2,3$, (see Lemma 5.3 at
the end of this subsection). Regarding $(5.30a)$ and summing it over
$i$ and $s$, we obtain

$$
\sum_{s,i}|\epsilon^{(2)}_{sp}|\le{Nc \sum_{|s|\le{n}}\{K_{4}{\bf
E}|H(s,p) D^{4}_{sp}[g^{0}_{1}G_{2}(p,s)]^{(1)} | +K_{2}{\bf
E}|H^{3}(s,p)D^{4}_{sp}[g^{0}_{1}G_{2}(p,s)]^{(2)} | }
$$
$$
+ {\bf E}|H^{5}(s,p)D^{4}_{sp}[g^{0}_{1}G_{2}(p,s)]^{(3)}| \} \eqno
(5.30b)
$$

 Hereafter, we use the notation:  for each pair $(p,s)$ and $\nu=1,2,3$, let
$H^{(\nu)}_{sp}=\hat{H}$ be the matrix defined by
  $$
  \hat{H}(r,i)= \left\{
\begin{array}{lll}
H(r,i) & \textrm{if} & (r,i)\neq(s,p) \\
\hat{H}(s,p) &  \textrm{if} & (r,i)=(s,p)
\end{array}\right.
$$
with $|\hat{H}(s,p)|\le{|H(s,p)|}$, the resolvent
$G^{(\nu)}_{sp}(z_{j})=\hat{G}_{j}$ and its normalized trace by
$\hat{g}_{j}$, $j=1,2$. To derive estimates of  $\epsilon^{(2)}$, we
need the following Lemma.
  \vskip0,5cm

 {\bf Lemma 5.2}. {\it If
$z\in\Lambda_{\eta}$ $(2.10)$ and for large enough $n$ and $b$, the
estimates
$$
D^{r}_{sp}\{ g^{0}_{n}(z) \} = O(\frac{1}{N}),  \quad r=1,2,.. \eqno
(5.31)
$$
and
$$
{\bf Var}(\hat{g}_{n}(z)) = O\{{\bf Var}(g_{n}(z)) +
\frac{1}{bN^{2}}\}, \eqno (5.32)
$$
are true in the limit  $N,b\longrightarrow \infty$. } \vskip0,5cm
 We prove Lemma 5.2 at the end of this section.
\vskip0,5cm
 Keeping in mind the expression $(5.30b)$ and using
$(5.28)$, $(5.31)$ and $(5.32)$, we obtain inequalities
$$
\sum_{|s|\le{n}}{\bf
E}|\hat{H}^{5}(s,p)D^{4}_{sp}\{\widehat{g_{1}^{0}G_{2}(p,s)}\}| \le{
\frac{\beta \hat{\mu}_{5}}{Nb^{3/2}}\sum_{|s|\le{n}}
\frac{\psi(\frac{s-p}{b})}{b}+\beta \sum_{|s|\le{n}}{\bf
E}\{|\hat{H}(s,p)|^{5}|\hat{g}^{0}_{1}|\} }
$$
$$
\le{  \frac{\beta \hat{\mu}_{5}}{Nb^{3/2}}\sum_{|s|\le{n}}
\frac{\psi(\frac{s-p}{b})}{b}+\beta \sqrt{{\bf
Var}\{\hat{g}_{1}\}}\sum_{|s|\le{n}} \sqrt{{\bf
E}\{|\hat{H}(s,p)|^{10}\}} }
$$
$$
\le{  \frac{\beta \hat{\mu}_{5}}{Nb^{3/2}}\sum_{|s|\le{n}}
\frac{\psi(\frac{s-p}{b})}{b}+\frac{\beta
\sqrt{\hat{\mu}_{10}}}{b^{3/2}} \sqrt{{\bf
Var}\{\hat{g}_{1}\}}\sum_{|s|\le{n}}\frac{\sqrt{\psi}(\frac{s-p}{b})}{b}
}
$$
$$
=O\{\frac{1}{b^{3/2}}[\frac{1}{N}+\sqrt{{\bf Var}\{g_{1}\}} ]\},
\eqno (5.33)
$$

$$
 \sum_{|s|\le{n}}K_{4}{\bf
E}|\hat{H}(s,p)D^{4}_{sp}\{\widehat{g_{1}^{0}G_{2}(p,s)}\}|
$$
$$
\le{\frac {3v^{2}\beta
\hat{\mu}_{1}}{Nb^{3/2}}\sum_{|s|\le{n}}\frac{\psi(\frac{s-p}{b})}{b}+
\frac{3v^{2}\beta \sqrt{\hat{\mu}_{2}}\sqrt{{\bf
Var}\{\hat{g}_{1}\}}}{b^{3/2}}\sum_{|s|\le{n}}\frac{\psi(\frac{s-p}{b})}{b}}
$$
$$
=O\{\frac{1}{b^{\frac{3}{2}}}[\frac{1}{N}+\sqrt{{\bf Var}\{g_{1}\}}
]\} \eqno (5.34)
$$
and
$$
 \sum_{|s|\le{n}}K_{2}{\bf
E}|\hat{H}^{3}(s,p)D^{4}_{sp}\{\widehat{g_{1}^{0}G_{2}(p,s)}\}|
$$
$$
\le{\frac {3v^{2}\beta
\hat{\mu}_{3}}{Nb^{3/2}}\sum_{|s|\le{n}}\frac{\psi(\frac{s-p}{b})}{b}+
\frac{3v^{2}\beta \sqrt{\hat{\mu}_{6}}\sqrt{{\bf
Var}\{\hat{g}_{1}\}}}{b^{3/2}}\sum_{|s|\le{n}}\frac{\psi(\frac{s-p}{b})}{b}}
$$
$$
=O\{\frac{1}{b^{\frac{3}{2}}}[\frac{1}{N}+\sqrt{{\bf Var}\{g_{1}\}}
]\} \eqno (5.35)
$$
with some constant $\beta$ and $\hat{\mu}_{r}=\sup_{i,p}{\bf
E}|a(i,p)|^{r}<\infty$, $r=1,..,10$. Then $(5.33)$, $(5.34)$ and
$(5.35)$ imply the following estimate

$$
|\sum_{s,i}\epsilon^{(2)}_{sp}| =
O\{N[\frac{1}{b^{\frac{3}{2}}}(\frac{1}{N} +\sqrt{{\bf
Var}(\hat{g}_{1})})]\}
$$
$$
=O[N(\frac{1}{Nb}+\frac{1}{b}\sqrt{{\bf Var}(g_{1})})]. \eqno (5.36)
$$

Returning to  $(5.23)$, using definition $N=2n+1>n$ and gathering
estimates $(5.24)$, $(5.25)$, $(5.26)$, $(5.27)$, $(5.29)$ and
$(5.36)$, we get $(5.21)$. \vskip1cm

 Now, let us prove $(5.22)$. Regarding $(5.23)$ and summing it over $p$ with
$G_{2}=\overline{G}_{1}=\overline{G}$ and using $(5.21)$, we obtain
estimate
$$
{\bf Var}\{g(z)\}\le{A(\frac{1}{nb}+\frac{1}{b}\sqrt{{\bf
Var}(g(z))})}. \eqno (5.37)
$$
where $A$ is a constant. Regarding this inequality, it is easy to
show that $ {\bf Var}\{g_{n}(z)\}=O(b^{-2}) $. Relation $(5.22)$ is
then proved and ends the proof of Theorem 5.2. $\diamond$ \vskip1cm

 {\it Proof of Lemma 5.2}. We start with $(5.31)$. Consider the
normalized trace $g_{n}(z)=g(z)=\frac{1}{N}\sum_{|k|\le{n}}G(k,k)$,
with $z\in\Lambda_{\eta}$. Regarding $(3.6)$, we obtain
$$
D^{1}_{sp}\{ g^{0}(z) \} = \frac{1}{N}\sum_{|k|\le{n}}D^{(1)}_{sp}\{
G(k,k) \}=-\frac{2}{N}\sum_{|k|\le{n}}G(k,s)G(k,p)
$$
$$
=-\frac{2}{N}G^{2}(s,p).
$$
Using $(3.7)$, we get
$$
|D^{1}_{sp}\{ g^{0}(z) \}|\le{\frac{2}{N|Imz|^{2}}}.
$$
Regarding the second derivative, one obtains
$$
D^{2}_{sp}\{ g^{0}(z) \} = \frac{2}{N}\{
2G^{2}(s,p)G(s,p)+G^{2}(p,p)G(s,s) +G^{2}(s,s)G(p,p)\},
$$
using $(3.7)$, we obtain $D^{2}_{sp}\{ g^{0}(z) \} =
O(\frac{1}{N})$. Then it is easy to show that $D^{r}_{sp}\{ g^{0}(z)
\}= O(\frac{1}{N})$, $r=1,2..$. Estimate $(5.31)$ is
proved.\vskip1cm

  Now let us prove $(5.32)$. Using the resolvent identity $(3.1)$ for
two hermitian matrices $H$ and $\hat{H}$ , we obtain
  $$
  \hat{g}=\hat{g}_{n}(z)=\frac{1}{N}\sum_{|k|\le{n}}\hat{G}(k,k)=\frac{1}{N}\sum_{|k|\le{n}}
  G(k,k)+\frac{1}{N}\sum_{k,r,i}\hat{G}(k,r)[\hat{H}-H](r,i)G(i,k)
  $$
$$
=g+\frac{1}{N} Tr(G\hat{G}\delta_{H}) \eqno (5.38)
$$
with
$$
[\delta_{H}](r,i)=[\hat{H}-H](r,i)=\left\{
\begin{array}{lll}
0 & \textrm{if} & (r,i)\neq(s,p), \\
\hat{H}(s,p)-H(s,p) &  \textrm{if} & (r,i)=(s,p).
\end{array}\right.
$$
Then
$$
Tr (G\hat{G}\delta_{H})=\sum_{r,i}
[G\hat{G}](r,i)[\delta_{H}](r,i)=[G\hat{G}](s,p)[\hat{H}(s,p)-H(s,p)].
\eqno (5.39)
$$

Let us return to relation $(5.38)$. It is easy to see that
$$
{\bf Var}(\hat{g})\le{ 2{\bf Var}(g) + \frac{2}{N^{2}}{\bf
Var}(Tr\{G\hat{G}\delta_{H}\})    }. \eqno (5.40)
$$
Then assuming that
$$
{\bf Var}(Tr\{G\hat{G}\delta_{H}\})=O(\frac{1}{b}) \eqno (5.41)
$$
holds, $(5.32)$ follows from $(5.40)$ and $(5.41)$. Let us prove
$(5.41)$. It is clear that $ |\hat{H}_{sp}-H_{sp}|\le{2|H_{sp}|}$.
Using $(3.7)$ and $(5.39)$, we obtain
$$
{\bf Var}(Tr\{G\hat{G}\delta_{H}\})\le{\frac{4}{|Imz|^{4}}{\bf
E}(|H(s,p)|+{\bf E}(|H(s,p)|)   )^{2}  }
$$
$$
\le{\frac{4}{b|Imz|^{4}}\{[{\bf E}(|a_{sp}|) (1 +
\psi(\frac{s-p}{b}))  ]^{2}\psi(\frac{s-p}{b})+{\bf
E}(|a_{sp}|)^{2}\psi(\frac{s-p}{b})^{2}(1-\psi(\frac{s-p}{b}))\} }
$$
$$
=O(\frac{1}{b}). \quad  \diamond
$$

Finally, we prove $(5.30)$ for the remainder $\epsilon^{(2)}$.
 \vskip1cm

 {\bf Lemma 5.3}. ({\bf The cumulant expansions method with q=3}).
 {\it  Let us consider the
family $\{X_{j}, \ j=1,..,m\}$ of independent real random variables
determined on the same probability space such that ${\bf
E}\{|X_{j}|^{5}\}<\infty $, $j=1,..,m$ and
$$
{\bf E}(X_{j})= {\bf E}(X_{j}^{3})=0, \quad j=1,..,m.
$$
If $F(t_{1},..,t_{m})$ is a complex-valued function acting on
vectors (with components as a real variables) in a $m-$dimensional
space such that its first $4$ derivatives are continuous and
bounded, then, for all $j$, we obtain $(3.4)$ with $q=3$ and the
remainder $\epsilon_{3}$ is given by
$$
\epsilon_{3}=-\frac{K_{4}}{3!}{\bf E}\{Xf_{j}^{(4)}(Y_{2})\}
-\frac{K_{2}}{3!}{\bf E}\{X^{3}f_{j}^{(4)}(Y_{1})\}
$$
$$
+\frac{1}{4!}{\bf E}\{X^{5}f_{j}^{(4)}(Y_{0})\}, \eqno (5.42)
$$
where $|Y_{\nu}|\le{|X_{j}|}$, $\nu=0,1,2$, $K_{r}$ is the r-th
cumulant of $X_{j}$ and $f_{j}$ is a complex-valued function of one
real variable such that
$$
f_{j}(Y)=F(X_{1},..,X_{j-1},Y,X_{j+1},..,X_{n})
$$
and $f_{j}^{(4)}$ is its $4-$th derivative.

 The cumulants can be expressed in terms of the moments. If $\mu_{r}={\bf
E}(X_{j}^{r})$ with $j=1,..,m$, one obtains: $K_{1}=\mu_{1}=0$,
$K_{2}=\mu_{2}$, $K_{3}=0$, $K_{4}=\mu_{4}-3\mu^{2}_{2}$.
 }
\vskip0,5cm

 {\it Proof of Lemma 5.3}. To simplify the notation we write $X=X_{j}$ for a
real random variable such that ${\bf
 E}\{|X|^{5}\}<\infty$ and ${\bf E}(X)={\bf E}(X^{3})=0$, and we denote
$f=f_{j}$.

 Now, consider a function $f$ of the class $C^{4}$ and apply Taylor's formula
to $f$, to obtain relation
$$
f(X)=\sum_{r=0}^{3}\frac{X^{r}}{r!}f^{(r)}(0)+\frac{X^{4}}{4!}f^{(4)}(Y_{0}),
 \quad |Y_{0}|\le{|X|}.
$$
Multiplying this equation by $X$ and taking the mathematical
expectation of both sides, we obtain equality
$$
{\bf E}\{Xf(X)\}={\bf E}(X^{2}){\bf E}(f^{'}(X))+{\bf
E}(\frac{X^{4}}{3!})f^{(3)}(0)+{\bf
E}\{\frac{X^{5}}{4!}f^{(5)}(Y_{0})\}. \eqno (5.43)
$$
Now we apply Taylor's formula to $f^{'}$, and write that
$$
f^{'}(X)=f^{'}(0)+Xf^{(2)}(0)+\frac{X^{2}}{2!}f^{(3)}(0)+
\frac{X^{3}}{3!}f^{(4)}(Y_{1}),
$$
with $|Y_{1}|\le{|X|}$. Taking the mathematical expectation of both
sides and multiplying by ${\bf E}(X^{2})$, we obtain relation
$$
{\bf E}(X^{2}){\bf E}\{f^{'}(0)\}={\bf E}(X^{2}){\bf
E}\{f^{'}(X)\}-\frac{[{\bf E}(X^{2})]^{2}}{2!}{\bf E}\{f^{(3)}(0)\}
$$
$$
-{\bf E}(X^{2}){\bf E}\{\frac{X^{3}}{3!}f^{(4)}(Y_{1})\}. \eqno
(5.44)
$$
Applying Taylor's formula to $f^{(3)}$, we see that
$$
{\bf E}\{f^{(3)}(0)\}={\bf E}\{f^{(3)}(X)\}-{\bf
E}\{Xf^{(4)}(Y_{2})\}
$$
with $|Y_{2}|\le{|X|}$. Then we can write that
$$
{\bf E}(X^{4}){\bf E}\{f^{(3)}(0)\}={\bf E}(X^{4}){\bf
E}\{f^{(3)}(X)\}- {\bf E}(X^{4}){\bf E}\{Xf^{(4)}(Y_{2})\} \eqno
(5.45)
$$
and that
$$
{\bf E}(X^{2})^{2}{\bf E}\{f^{(3)}(0)\}={\bf E}(X^{2})^{2}{\bf
E}\{f^{(3)}(X)\}- {\bf E}(X^{2})^{2}{\bf E}\{Xf^{(4)}(Y_{2})\}.
\eqno (5.46)
$$
Now, substitute $(5.46)$ in relation $(5.44)$ and use $(5.45)$, to
rewrite $(5.43)$ as
$$
{\bf E}\{Xf(X)\}={\bf E}(X^{2}){\bf
E}(f^{'}(X))+\frac{K_{4}}{3!}{\bf E}\{f^{(3)}(X)\}
$$
$$
+\frac{1}{4!}{\bf E}\{X^{5}f^{(4)}(Y_{0})\}-\frac{K_{4}}{3!}{\bf
E}\{Xf^{(4)}(Y_{2})\}-\frac{K_{2}}{3!}{\bf E}\{X^{3}f^{(4)}(Y_{1})\}
$$
with $K_{4}=[{\bf E}(X^{4})-3{\bf E}(X^{2})^{2}]$ and $K_{2}={\bf
E}(X^{2})$. This gives relation $(5.42)$. $\diamond$

\vskip0,5cm

 {\bf Acknowledgements.}
The author is grateful to Prof. Dr. O. Khorunzhy at University of Versailles (France), where present paper was completed, who proposed use to study the problems described in this paper.

\end{document}